\documentclass{article}

\usepackage[T1]{fontenc}
\usepackage{lmodern}

\usepackage{authblk}
\usepackage[comma]{natbib}
\usepackage{graphicx}
\usepackage{amssymb}

\usepackage{amsmath,amsfonts,bm}

\def\eqref#1{equation~\ref{#1}}

\def\1{\bm{1}}

\def\vb{{\bm{b}}}
\def\vc{{\bm{c}}}

\def\vh{{\bm{h}}}

\def\vl{{\bm{l}}}

\def\vu{{\bm{u}}}

\def\vx{{\bm{x}}}
\def\vy{{\bm{y}}}
\def\vz{{\bm{z}}}

\def\mA{{\bm{A}}}

\def\mW{{\bm{W}}}

\DeclareMathAlphabet{\mathsfit}{\encodingdefault}{\sfdefault}{m}{sl}
\SetMathAlphabet{\mathsfit}{bold}{\encodingdefault}{\sfdefault}{bx}{n}

\def\sR{{\mathbb{R}}}

\usepackage{tikz}
\usetikzlibrary{positioning,calc}
\usepackage{enumerate}
\usepackage{multicol}
\usepackage{hyperref}
\usepackage{doi}
\usepackage{multirow}
\usepackage{placeins}

\newcommand{\beginFEappendix}{\appendix}
\newcommand{\finishFEappendix}{ }

\title{Constraint Learning by Optimizing over Trained Neural Networks: A Pragmatic Tutorial}

\author{Jiaxiao Fang}
\author{Thiago Serra}

\affil{University of Iowa}

\date{October 2026}

\begin{document}

\maketitle

\begin{abstract}
\noindent
Constraint learning is an emerging paradigm of modeling optimization problems with the help of machine learning. 
In this paradigm, trained machine learning models are embedded into the optimization model.
For example, we can learn constraints empirically from available data, or we can avoid using a nonlinear objective function by learning a simpler approximation for it. 
Neural networks become a compelling option if the machine learning models must capture more elaborate and nuanced relationships. 
However, optimization models embedding neural networks may become too large or too slow to solve. 
In this tutorial, we present practical considerations for an effective and efficient use of constraint learning with neural networks 
and validate it with a case study highlighting the benefit of smaller and sparser networks. 

\end{abstract}

\section{Introduction}\label{sec:intro}

We may not fully understand how neural networks arrive at a conclusion, 
but we know that they can handle complex tasks better than more intuitive counterparts, 
such as decision trees~\citep{baldi2014deephep,lecun2015deep}. 
That fits into a broader narrative outlined in \citeauthor{sutton2019bitter}'s manifesto The Bitter Lesson 
about how general-purpose methods backed by massive computation led to the most effective advances in AI, 
from playing chess and Go at competitive levels to achieving reliable results in speech recognition and computer vision~\citep{sutton2019bitter}. 
And this narrative begets the question: 
\emph{what if some parts of a mathematical optimization model were learned with neural networks?} 
In machine learning, 
the cost of producing a model (\emph{training}) is generally orders of magnitude above the cost of using it (\emph{inference}). 
However, 
if neural networks are embedded into a mathematical optimization model  
and possibly queried several times while solving the optimization model, 
then the inference cost adds up to another order of magnitude. 
This tutorial is about \emph{when} to use neural networks in optimization models---and \emph{how} to circumvent their computational cost.

We define as \emph{constraint learning} the use of neural networks, or other machine learning models, 
as a proxy for some parts of a mathematical optimization model. 
However, that term is also used more generally for the automatic process of learning the model of a problem in Constraint Programming~(CP)\footnote{CP is another technique for solving discrete optimization problems. For more information, see \nameref{ap:about_cp} (in Appendix).}~\citep{kolb2016learning,lombardi2018survey}.  
As observed by \cite{deraedt2018learningconstraints}, 
constraint learning stands for an \emph{inductive} (probabilistic) approach of learning constraints, 
in contrast to the \emph{deductive} (provable) approach used for inferring  valid constraints (cuts). 
When machine learning is involved, 
this probabilistic approach translates into \emph{approximating} through data a complex or unknown constraint or term of the objective function. 
The special case of using neural networks has been first explored in CP models by~\cite{bartolini2011neuron} 
and in mathematical optimization models by \cite{say2017nonlinear} and \cite{cheng2017maximum}, 
to which many other studies followed in both communities. 

Curiously, 
the term constraint learning would only find its way into mathematical optimization six years later, but would soon become widely used as a synonym for embedding neural networks in optimization models.\footnote{For a list of references, see \nameref{ap:adopting_cl} (in Appendix).}
However, other terms have been used to denote special or general cases of constraint learning.
While \cite{lombardi2017edml} suggested denoting the use of general-purpose machine learning models for constraint learning as \emph{empirical model learning},  
the practice of using simpler models to approximate complex or unknown relationships in engineering has been generally denoted as 
\emph{surrogate modeling}~\citep{bhosekar2018review,grimstad2019surrogate}.

Naming controversies aside, 
there are some important  considerations for achieving 
a successful use of constraining learning 
with neural networks. 
This tutorial covers them in the following order:
\begin{itemize}
\item We discuss why neural networks can be helpful; what type of neural network to use and how to incorporate them in mathematical optimization models; 
and where they have been applied in \nameref{sec:why_what_where} This is our summary of positives. 
\item We briefly discuss when neural networks exceed the modeling needs and we present alternatives for such cases in \nameref{sec:not} This is our warning. 
\item Most importantly, we discuss what technical choices lead to a more effective use by illustrating those choices with experiments on a synthetic dataset in \nameref{sec:how} This is our practical advice, which we support with experimental results. 
\item Finally, we briefly discuss advanced topics and related areas beyond the scope of this tutorial in \nameref{sec:etc}. 
\end{itemize}

\section{Why, What, and Where to Use Neural Networks in Modeling?}\label{sec:why_what_where}

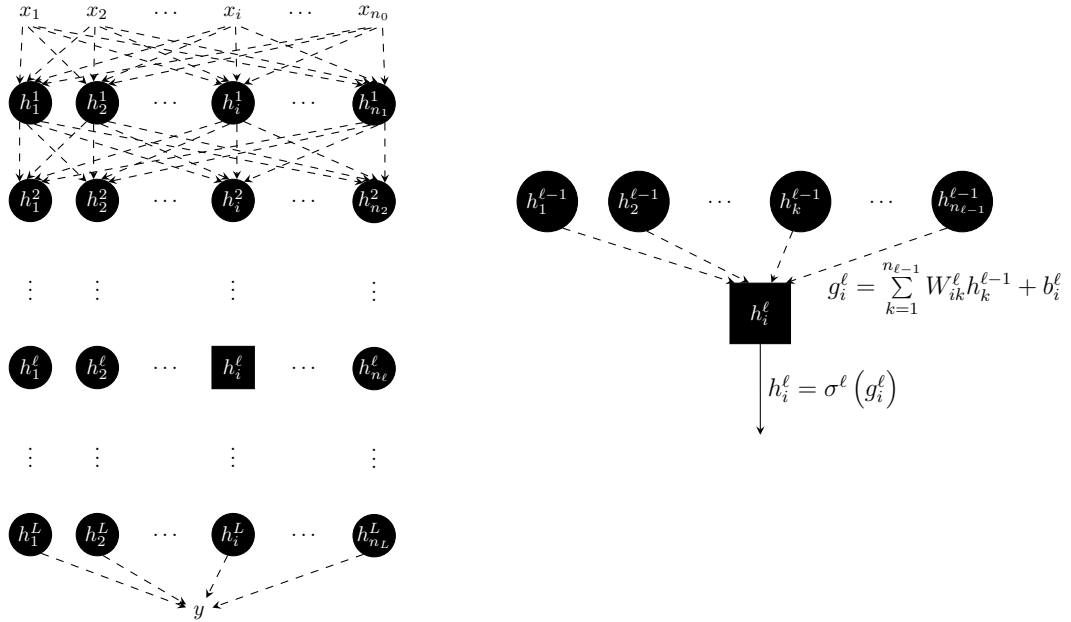
\begin{figure}[h!]
    \centering
    \newcommand{\nodesep}{0.5cm}
\newcommand{\layersep}{0.9cm}
\newcommand{\circlesize}{7mm}

\begin{tikzpicture}[>=stealth,scale=0.8,transform shape]
\tikzset{
  neuron/.style={draw,circle,fill=black,minimum size=\circlesize,inner sep=0pt},
  conn/.style={->,dashed,align=left},
  conn2/.style={->,align=left}
}

\node (x1) {$x_1$};
\node (x2) [right=\nodesep of x1] {$x_2$};
\node (e1) [right=\nodesep of x2] {$\ldots$};
\node (xj) [right=\nodesep of e1] {$x_i$};
\node (e2) [right=\nodesep of xj] {$\ldots$};
\node (xn) [right=\nodesep of e2] {$x_{n_0}$};

\node[neuron] (h11) [below=\layersep of x1] {};
\node[neuron] (h12) [below=\layersep of x2] {};
\node[neuron] (h1j) [below=\layersep of xj] {};
\node[neuron] (h1n) [below=\layersep of xn] {};

\path (h12.center) -- (h1j.center) node[midway] (he11) {$\ldots$};
\path (h1j.center) -- (h1n.center) node[midway] (he12) {$\ldots$};

\node[text=white] at (h11) {$h^1_1$};
\node[text=white] at (h12) {$h^1_2$};
\node[text=white] at (h1j) {$h^1_i$};
\node[text=white] at (h1n) {$h^1_{n_1}$};

\node[neuron] (h21) [below=\layersep of h11] {};
\node[neuron] (h22) [below=\layersep of h12] {};
\node[neuron] (h2j) [below=\layersep of h1j] {};
\node[neuron] (h2n) [below=\layersep of h1n] {};

\path (h22.center) -- (h2j.center) node[midway] (he21) {$\ldots$};
\path (h2j.center) -- (h2n.center) node[midway] (he22) {$\ldots$};

\node[text=white] at (h21) {$h^2_1$};
\node[text=white] at (h22) {$h^2_2$};
\node[text=white] at (h2j) {$h^2_i$};
\node[text=white] at (h2n) {$h^2_{n_2}$};

\node (hvi1) [below=\layersep of h21] {};
\node (hvi2) [below=\layersep of h22] {};
\node (hvij) [below=\layersep of h2j] {};
\node (hvin) [below=\layersep of h2n] {};

\node at (hvi1) {$\vdots$};
\node at (hvi2) {$\vdots$};
\node at (hvij) {$\vdots$};
\node at (hvin) {$\vdots$};

\node[neuron] (hi1) [below=\layersep of hvi1] {};
\node[neuron] (hi2) [below=\layersep of hvi2] {};
\node[draw,rectangle,fill=black,minimum size=\circlesize,inner sep=0pt] (hij) [below=\layersep of hvij] {};
\node[neuron] (hin) [below=\layersep of hvin] {};

\path (hi2.center) -- (hij.center) node[midway] (hei1) {$\ldots$};
\path (hij.center) -- (hin.center) node[midway] (hei2) {$\ldots$};

\node[text=white] at (hi1) {$h^\ell_1$};
\node[text=white] at (hi2) {$h^\ell_2$};
\node[text=white] at (hij) {$h^\ell_i$};
\node[text=white] at (hin) {$h^\ell_{n_\ell}$};

\node (hv2i1) [below=\layersep of hi1] {};
\node (hv2i2) [below=\layersep of hi2] {};
\node (hv2ij) [below=\layersep of hij] {};
\node (hv2in) [below=\layersep of hin] {};

\node at (hv2i1) {$\vdots$};
\node at (hv2i2) {$\vdots$};
\node at (hv2ij) {$\vdots$};
\node at (hv2in) {$\vdots$};

\node[neuron] (hL1) [below=\layersep of hv2i1] {};
\node[neuron] (hL2) [below=\layersep of hv2i2] {};
\node[neuron] (hLj) [below=\layersep of hv2ij] {};
\node[neuron] (hLn) [below=\layersep of hv2in] {};

\path (hL2.center) -- (hLj.center) node[midway] (heL1) {$\ldots$};
\path (hLj.center) -- (hLn.center) node[midway] (heL2) {$\ldots$};

\node[text=white] at (hL1) {$h^L_1$};
\node[text=white] at (hL2) {$h^L_2$};
\node[text=white] at (hLj) {$h^L_i$};
\node[text=white] at (hLn) {$h^L_{n_L}$};

\node[inner sep=0pt] (y2) [below=\layersep of hL2] {};
\node[inner sep=0pt] (yj) [below=\layersep of hLj] {};
\node (y) at ($(y2.center)!0.75!(yj.center)$) {$y$};

\draw[conn] (x1.240) -- (h11.120);
\draw[conn] (x1.260) -- (h12.120);
\draw[conn] (x1.280) -- (h1j.120);
\draw[conn] (x1.300) -- (h1n.120);

\draw[conn] (x2.240) -- (h11.100);
\draw[conn] (x2.260) -- (h12.100);
\draw[conn] (x2.280) -- (h1j.100);
\draw[conn] (x2.300) -- (h1n.100);

\draw[conn] (xj.240) -- (h11.80);
\draw[conn] (xj.260) -- (h12.80);
\draw[conn] (xj.280) -- (h1j.80);
\draw[conn] (xj.300) -- (h1n.80);

\draw[conn] (xn.240) -- (h11.60);
\draw[conn] (xn.260) -- (h12.60);
\draw[conn] (xn.280) -- (h1j.60);
\draw[conn] (xn.300) -- (h1n.60);

\draw[conn] (h11.240) -- (h21.120);
\draw[conn] (h11.260) -- (h22.120);
\draw[conn] (h11.280) -- (h2j.120);
\draw[conn] (h11.300) -- (h2n.120);

\draw[conn] (h12.240) -- (h21.100);
\draw[conn] (h12.260) -- (h22.100);
\draw[conn] (h12.280) -- (h2j.100);
\draw[conn] (h12.300) -- (h2n.100);

\draw[conn] (h1j.240) -- (h21.80);
\draw[conn] (h1j.260) -- (h22.80);
\draw[conn] (h1j.280) -- (h2j.80);
\draw[conn] (h1j.300) -- (h2n.80);

\draw[conn] (h1n.240) -- (h21.60);
\draw[conn] (h1n.260) -- (h22.60);
\draw[conn] (h1n.280) -- (h2j.60);
\draw[conn] (h1n.300) -- (h2n.60);

\draw[conn] (hL1.300) -- (y);
\draw[conn] (hL2.285) -- (y);
\draw[conn] (hLj.255) -- (y);
\draw[conn] (hLn.240) -- (y);

\renewcommand{\circlesize}{10mm}
\renewcommand{\layersep}{1.5cm}

\node[neuron] (hi11) [right=\nodesep*4 of h2n] {};
\node[neuron] (hi12) [right=\nodesep of hi11] {};
\node (hi12e21) [right=\nodesep of hi12] {$\ldots$};
\node[neuron] (hi1j) [right=\nodesep of hi12e21] {};
\node (hi12e22) [right=\nodesep of hi1j] {$\ldots$};
\node[neuron] (hi1n) [right=\nodesep of hi12e22] {};

\node[text=white] at (hi11) {$h^{\ell-1}_1$};
\node[text=white] at (hi12) {$h^{\ell-1}_2$};
\node[text=white] at (hi1j) {$h^{\ell-1}_k$};
\node[text=white] at (hi1n) {$h^{\ell-1}_{n_{\ell-1}}$};

\node[inner sep=0pt] (hij2A) [below=\layersep of hi12e21] {};
\node[inner sep=0pt] (hij2B) [below=\layersep of hi1j] {};
\node[draw,rectangle,fill=black,minimum size=\circlesize,inner sep=0pt] (hij2) at ($(hij2A.center)!0.5!(hij2B.center)$) {};

\node[text=white] at (hij2) {$h^{\ell}_i$};

\node[inner sep=0pt] (invisible) [below=\layersep of hij2] {};

\draw[conn] (hi11.300) -- (hij2.130);
\draw[conn] (hi12.285) -- (hij2.110);
\draw[conn] (hi1j.255) -- (hij2.70);
\draw[conn] (hi1n.240) -- node[midway,right] {\\  \vspace{-1cm} \hspace{-1cm} \large $g_i^\ell =  \sum\limits_{k=1}^{n_{\ell-1}} W^\ell_{i k} h^{\ell-1}_k + b^\ell_i$}  (hij2.50);

\draw[conn2] (hij2.270) -- node[midway,right] {\large $h^\ell_i = \sigma^\ell \left( g^\ell_i \right)$} (invisible.90);

\end{tikzpicture}
    \caption{Left: Structure of a feedforward neural network with input $\vx \in \sR^{n_0}$, hidden layers $\ell=1$ to $L$ with $n_\ell$ neurons each, and output $y \in \sR$. Right: Formula for output $h^\ell_i$ of neuron $i$ in layer $\ell$ from  activation function $\sigma^\ell : \sR \rightarrow \sR$ and preactivation output $g^\ell_i$, which uses previous layer output $\vh^{\ell-1}$, neuron weights $\mW^\ell_{i *}$, and bias $b^\ell_i$. For details, see \nameref{ap:nn_structure_map} (in Appendix).}
    \label{fig:nn_figure}
\end{figure}

For our purposes, it is convenient to think of a trained neural network as a function $\hat{f} : \sR^{n_0} \rightarrow \sR$ 
used to approximate an unknown function $f$, 
and so we may refer to $\hat{f}$ as a \emph{proxy function} with respect to $f$. 
Moreover, it is sufficient to consider a feedforward neural network, such as illustrated in Figure~\ref{fig:nn_figure}. 
There are certainly more sophisticated architectures that could be used for a neural network, but we will justify along this section why we do not recommend starting with them.

In the case of an optimization problem over decision variables in $\vx$ that is otherwise linear  
and with decision variables that are continuous ($\vx \in \mathbb{R}^n$), integer ($\vx \in \mathbb{Z}^n$), or mixed-integer ($\vx \in \mathbb{Z}^p \times \mathbb{R}^{n-p}$), 
a neural network can model a constraint or term of the objective function as follows:
\begin{multicols}{2}
\noindent
\begin{align}
\begin{aligned}\label{eq:f_cons}
    \max & ~ & \vc^T \vx \\
    \text{s.t.} & ~ & \mA \vx \leq \vb \\
    && \hat{f}(\vx) \leq  d 
\end{aligned}
\end{align}
\begin{align}
\begin{aligned}\label{eq:f_obj}
    \max & ~ & \vc^T \vx + \hat{f}(\vx) \\
    \text{s.t.} & ~ & \mA \vx \leq \vb 
\end{aligned}
\end{align}
\end{multicols}
\noindent
For example, we can think of $\hat{f}$ as an estimate of resource consumption capped at some value $d$ in (\ref{eq:f_cons}) or of profit in (\ref{eq:f_obj}). Note that $\hat{f}$ does not necessarily need to involve all decision variables in $\vx$. 

\subsection{Why Use a Neural Network?}

\emph{There are two main types of applications for constraint learning:}
\begin{enumerate}[(i)]
    \item \emph{When there is no formula mapping $x$ to $y$, but we can estimate it from the available data.} 
    \item \emph{When the formula mapping $x$ to $y$ is known, but it cannot be easily handled by a solver if stated in the optimization model. In this case, we can generate data and treat it as in (i).}
\end{enumerate}

Hence, we resort to machine learning to formulate the part of a mathematical optimization model that is either (i) \emph{unknown} or (ii) \emph{intractable}. 
Working with an unknown relationship based on the available data, 
such as in forecasting, 
is a common application of machine learning. 
Working with parts of a model that make an optimization problem harder to solve, 
such as with certain nonlinear functions or probability distributions, 
is an area of growing application of machine learning. 
We survey the constraint learning literature related to both in \nameref{sec:where}

What makes neural networks particularly appealing in this context is that they are known to be \emph{universal function approximators}, meaning that they can approximate most functions to an chosen precision $\varepsilon$. For example, suppose that we want to approximate a function $f$ with a neural network that produces a  function $\hat{f}$ after being trained. 
For any valid $\vx$, approximating to such a precision means that $|f(\vx) - \hat{f}(\vx)| \leq \varepsilon$. 
In other words, 
it is \emph{theoretically} possible that the output of the neural network will not differ from $f$ by more than $\varepsilon$. 
Finding such a neural network has been shown to be possible if there is a minimum number of neurons per layer (\emph{width}) and a sufficient number of layers (\emph{depth}); or, conversely, a minimum depth and a sufficient width instead~\citep{huchette2026polyhedral}. 
However, 
we may end up working with neural networks that might not be as wide or deep; and we may depend on training algorithms that have limited theoretical guarantees.  
Nevertheless, the successful use of neural networks across a wide variety of domains is evidence that we may still obtain very satisfactory results \emph{in practice}.

\begin{figure}[h!]
    \centering
    \includegraphics[width=0.9\linewidth]{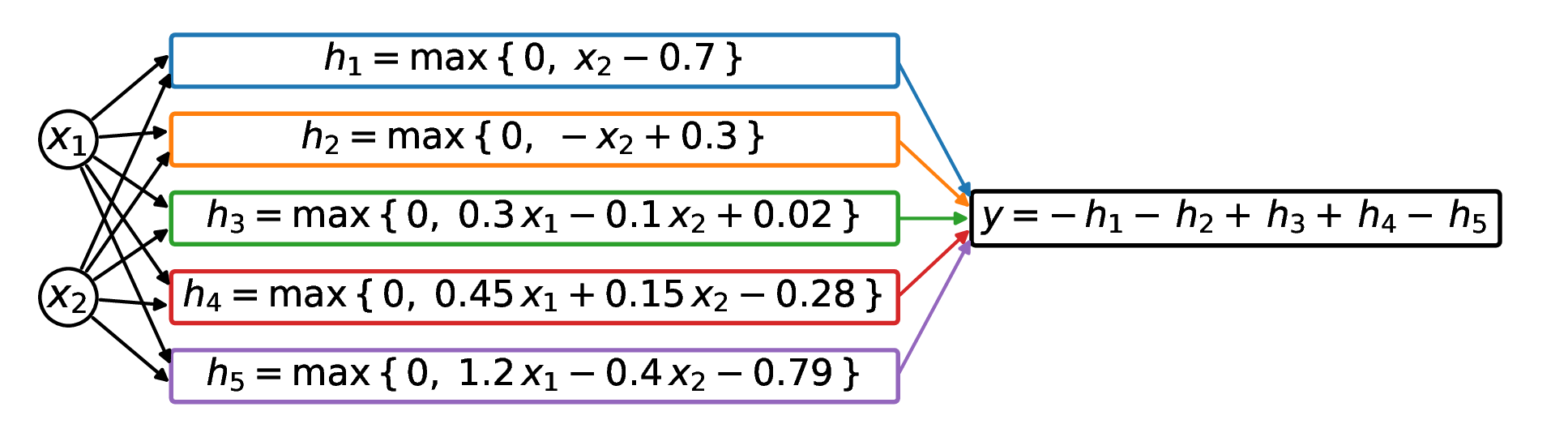}
    \includegraphics[scale=0.45]{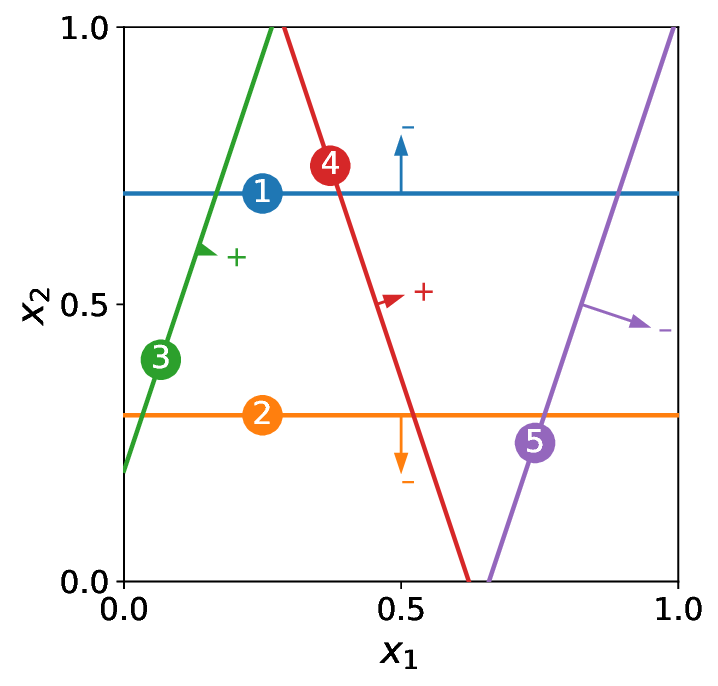}
    \includegraphics[scale=0.45]{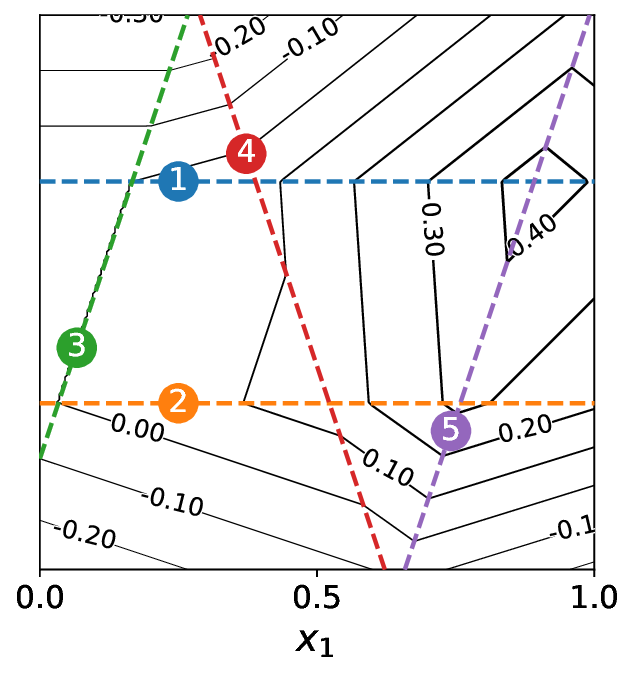}
    \caption{
    \textbf{Top:} Visual description of neural network used as example. 
    \textbf{Bottom left:} Lines partitioning the space based on what inputs produce a positive output for each neuron, with the arrow pointing to the positive side, the length of the arrow proportional to the magnitude of the weights, and the arrow label denoting the influence on $y$.
    \textbf{Bottom right:} Contour plots of $y$ over the lines associated with neural activations. 
    The example and plots are adapted from \cite{tong2026gradient}.
    }
    \label{fig:pwa}
\end{figure}

By considering neural networks as  approximators in constraint learning, 
we are hoping to find a compromise to make it easier to solve optimization problems. 
In fact, 
approximations are already widely used in nonlinear and stochastic optimization. 
In nonlinear optimization, 
piecewise-linear functions are often used to approximate nonlinear functions~\citep{mccormick1976under,mangasarian2005under,vielma2010pwl,misener2010multi,huchette2023pwl}, 
and we can indeed train neural networks to produce piecewise-linear functions. 
We illustrate in Figure~\ref{fig:pwa} an example of a neural network on $\vx \in [0,1]^2$ (top) whose output $y$ is piecewise-linear and depicted by contour plots (bottom right). 
We leave other parts of this figure, as well as the case of modeling probability distributions with neural networks, to be discussed in \nameref{sec:what_type}

\subsection{What Type of Neural Network Should I Use?}\label{sec:what_type}

\emph{In short, 
we recommend a feedforward neural network in which the activation function is the Rectified Linear Unit~(ReLU): $\sigma^\ell(u) = \max\{0, u\}$.} 
Now we proceed with the longer answer. 

\paragraph{Why is ReLU a Good Choice?} 

When we embed a neural network into an optimization model, 
we are effectively simulating the neural network within that model. 
However, 
unless the activation function $\sigma^\ell$ 
is a nonlinear function, 
the neural network is no more expressive than a linear regression.\footnote{
For an explanation, see \nameref{app:nonlinear_act} (in Appendix).}  
But we should be careful with our choice for $\sigma^\ell$ 
because nonlinear optimization models can be considerably more difficult to solve. 
As mentioned before, 
piecewise-linear approximations are commonly used to approach nonlinear optimization problems. 
That is because we can formulate piecewise-linear functions using only binary variables and linear constraints~\citep{vielma2010pwl}, 
as we illustrate after this discussion with formulation~(\ref{eq:relu_mip}). 
Hence, 
we can keep using only linear constraints and therefore a formulation based on (\ref{eq:f_cons}) or (\ref{eq:f_obj}) yields 
a Mixed-Integer Linear Programming~(MILP) model, 
which is more tractable than its nonlinear counterpart. 

The ReLU is convenient because it is a popular and relatively simple piecewise-linear function, 
consisting of the maximum between zero and the preactivation output of the neuron. 
ReLU started being adopted in deep neural networks by  \cite{nair2010relu} and \cite{glorot2011deep} and in few years became the most commonly used activation function~\citep{lecun2015deep,ramachandran2018searching}, only losing space with the emergence of the architectures currently used for generative artificial intelligence in transformers~\citep{vaswani2017attention} and diffusion models~\citep{ho2020ddpm}. 
In machine learning folklore, 
the ReLU is ``nonlinear enough, but no more than necessary.'' 

By ``nonlinear enough'', 
we mean that with the ReLU we can already model complex and nuanced relationships. 
Let us turn back to Figure~\ref{fig:pwa} to understand this in more detail. 
At the \emph{top} of Figure~\ref{fig:pwa}, we have a hidden layer with 5 neurons, 
each with its own parameters influencing its output. 
For example, the output of the first neuron is $h^1_1 = \max\{0, x_2 - 0.7\}$, 
implying that its output is positive ($h^1_1 > 0$) only when $x_2 - 0.7 > 0$ (i.e., $x_2 > 0.7$). 
In other words, 
the first neuron has no influence on the neural network output when $x_2 \leq 0.7$. 
At the \emph{bottom left} of Figure~\ref{fig:pwa}, 
this translates into the horizontal blue line at $x_2 = 0.7$, 
which has an arrow pointing upwards to indicate the side where the neuron influences the network output.  
When we increase the value of $x_2$ beyond 0.7, 
an increment of $\delta$ to $x_2$ implies an increment of $\delta$ to $h^1_1$. 
Because of the negative influence of $h^1_1$ on the network output $y$, 
an increment of $\delta$ to $x_2$ implies a decrement of $\delta$ to $y$.  
At the \emph{bottom right} of Figure~\ref{fig:pwa}, 
we see that the contour plots above $x_2 = 0.7$ 
are more horizontal and that the network output goes down as we increase $x_2$. 
This reflects the influence of the first neuron on the network output. 

Considering the influence of all five neurons in Figure~\ref{fig:pwa}, 
we have a piecewise-linear function with 12 pieces in $\vx \in [0,1]^2$. 
Each one of those pieces, which is known as a \emph{linear region}\label{pag:linear_region} in machine learning, 
corresponds to a different subset of neurons influencing the network output. 
Hence, 
we may think that every extra neuron doubles the number of linear regions. 
That exponential growth may indeed happen at first, but it slows down when a layer has more neurons than inputs and becomes polynomial for a fixed number of layers, 
as we will discuss later. 
Still, the number of linear regions may keep growing exponentially with the number of layers~\citep{huchette2026polyhedral}. 
This expressive power makes ReLU networks very appealing for modeling and approximation. 

\paragraph{How to Embed a ReLU in a Mathematical Optimization Model?} 
When we consider optimization and neural networks, 
there is an interesting contrast between training and constraint learning: 
\begin{itemize}
\item In training, the \emph{inputs} are fixed and we change the \emph{parameters} (sets of weights and bias) to obtain values for the preactivation output and output that should match the training set.
\item In constraint learning, the \emph{parameters} are fixed and we change the \emph{inputs} to obtain values for the preactivation output and output that would produce solutions that are feasible and better. 
\end{itemize}

We believe that it is important for the reader to understand how a neural network gets embedded into an optimization model, 
especially for discussing what design choices are more effective afterwards. 
But later the reader may use frameworks that automate this embedding, 
such as MeLOn~\citep{schweidtmann2020melon}, JANOS~\citep{bergman2022janos}, OMLT~\citep{ceccon2022omlt}, OCL~\citep{fajemisin2024ocl}, MathOptAI.jl~\citep{parker2024formulations}, 
Gurobi Machine Learning~\citep{gurobi2025ml}, and PySCIPOpt-ML~\citep{turner2025pyscipoptml}.

\subparagraph{Modeling Each Neuron} 
We will describe two equivalent formulations for modeling each neuron, 
as presented in \cite{serra2018bounding} and \cite{fischetti2018dnnmilp}. 
The first uses only algebraic expressions, which makes it implementable in any MILP solver. 
The second also uses logic expressions, which can be very convenient if they are supported by the solver. 

For every neuron $i$ in every layer $\ell$, 
the row of weights $\mW^\ell_{i *}$ and the bias $b^\ell_i$ of the neuron are constant; 
whereas the input $\vh^{\ell-1}$, preactivation output $g^\ell_i$, and output $h^\ell_i$ are decision variables. 
Hence, 
we formulate the algebraic mapping from the input $\vh^{\ell-1}$ to the output $h^\ell_i$ 
as follows:  
\begin{subequations}
\begin{align}
\sum_{k=1}^{n_{\ell-1}} \mW^\ell_{i k} h^{\ell-1}_k + b^\ell_i = g^\ell_i \label{eq:W} \\
g^\ell_i = h^\ell_i - \bar{h}^\ell_i \label{eq:g} \\
h^\ell_i \leq H^\ell_i z^\ell_i \label{eq:z} \\
\bar{h}^\ell_i \leq \bar{H}^\ell_i (1- z^\ell_i) \label{eq:znot} \\
h^\ell_i \geq 0 \label{eq:h} \\
\bar{h}^\ell_i \geq 0 \label{eq:hbar} \\
z^\ell_i \in \{0, 1\} \label{eq:bin}
\end{align}
\label{eq:relu_mip}
\end{subequations} 
With constraint (\ref{eq:W}), 
we model how the output $h_k^{\ell-1}$ of every neuron $k$ in layer $\ell-1$ contributes to the preactivation output $g_i^\ell$ of neuron $i$ in layer $\ell$. 
With the remaining constraints, (\ref{eq:g}) to (\ref{eq:bin}), 
we model the relationship between that preactivation output $g_i^\ell$ and the neuron output $h_i^\ell$. 
The central element in those constraints is the binary decision variable $z^\ell_i$ modeling if the neuron is \emph{active}:
\begin{itemize}
    \item If we set $z^\ell_i = 1$, then constraints (\ref{eq:g}) to (\ref{eq:bin}) simplify as follows:
    \begin{align}
    \small
    \left\{
    \begin{array}{c}
    g^\ell_i = h^\ell_i - \bar{h}^\ell_i \\
h^\ell_i \leq H^\ell_i \\
\bar{h}^\ell_i \leq 0 \\
h^\ell_i \geq 0  \\
\bar{h}^\ell_i \geq 0  \\
(z^\ell_i = 1)
    \end{array}
    \right\}
    \Rightarrow
    \left\{
    \begin{array}{c}
    g^\ell_i = h^\ell_i - \bar{h}^\ell_i \\
0 \leq h^\ell_i \leq H^\ell_i \\
\bar{h}^\ell_i = 0 
    \end{array}
    \right\}
    \Rightarrow
    \left\{
    \begin{array}{c}
    g^\ell_i = h^\ell_i \\
0 \leq h^\ell_i \leq H^\ell_i \\
(\therefore g^\ell_i \geq 0)
    \end{array}
    \right\}
    \end{align}
In the first step, we reduce the terms that are dependent on $z_i^\ell$, turning $H_i^\ell z_i^\ell$ into $H_i^\ell$ and $\bar{H}_i^\ell (1 - z_i^\ell)$ into $0$. 
In the second step, we combine the constraints bounding $h_i^\ell$ and $\bar{h}_i^\ell$ respectively into $0 \leq h^\ell_i \leq H^\ell_i$ and $\bar{h}_i^\ell = 0$. 
Finally, 
we state that $g_i^\ell = h_i^\ell$ because $\bar{h}_i^\ell=0$, 
which implies that $g_i^\ell \geq 0$.

    \item Otherwise, if we set $z^\ell_i = 0$, then constraints (\ref{eq:g}) to (\ref{eq:bin}) simplify as follows:
    \begin{align}
    \small
    \left\{
    \begin{array}{c}
    g^\ell_i = h^\ell_i - \bar{h}^\ell_i \\
h^\ell_i \leq 0 \\
\bar{h}^\ell_i \leq \bar{H}^\ell_i \\
h^\ell_i \geq 0  \\
\bar{h}^\ell_i \geq 0  \\
z^\ell_i = 1
    \end{array}
    \right\}
    \Rightarrow
    \left\{
    \begin{array}{c}
    g^\ell_i = h^\ell_i - \bar{h}^\ell_i \\
0 \leq \bar{h}^\ell_i \leq \bar{H}^\ell_i \\
h^\ell_i = 0 
    \end{array}
    \right\}
    \Rightarrow
    \left\{
    \begin{array}{c}
    g^\ell_i = - \bar{h}^\ell_i \\
0 \leq \bar{h}^\ell_i \leq \bar{H}^\ell_i \\
(\therefore g^\ell_i \leq 0)
    \end{array}
    \right\}
    \end{align}
\end{itemize}
In the first step, we reduce the terms that are dependent on $z_i^\ell$, turning $H_i^\ell z_i^\ell$ into $0$ and $\bar{H}_i^\ell (1 - z_i^\ell)$ into $\bar{H}_i^\ell$. 
In the second step, we combine the constraints bounding $h_i^\ell$ and $\bar{h}_i^\ell$ respectively into $h_i^\ell = 0$ and $0 \leq \bar{h}^\ell_i \leq \bar{H}^\ell_i$. 
Finally, 
we state that $g_i^\ell = - \bar{h}_i^\ell$ because $h_i^\ell=0$, 
which implies that $g_i^\ell \leq 0$.

In terms of the ReLU function $h^\ell_i = \max\{0, g^\ell_i\}$, 
$z^\ell_i = 1$ corresponds to the case in which the neuron is active and $h^\ell_i = g^\ell_i$; 
whereas $z^\ell_i=0$ corresponds to the case in which the neuron is inactive and $h^\ell_i = 0$.\footnote{A neuron is active when $g^\ell_i>0$ and inactive when $g^\ell_i<0$. But what if $g^\ell_i=0$? For our purposes, since $h^\ell_i = g^\ell_i = 0$, it can be either.} 
However, this formulation is only valid if, for any possible input $\vh^{\ell-1}$, we assure that $g_i^\ell \leq H_i^\ell$ in the first case (since $h_i^\ell \leq H_i^\ell$) and $g_i^\ell \geq - \bar{H}_i^\ell$ in the second case (since $\bar{h}_i^\ell \leq H_i^\ell$). 
Hence, the values of $H_i^\ell$ and $\bar{H}_i^\ell$ should be large enough. 
Those are the big--M constants used in the big--M constraints (\ref{eq:z}) and (\ref{eq:znot}): 
when $z_i^\ell = 1$, it is as if constraint (\ref{eq:z}) did not exist; 
and likewise for constraint (\ref{eq:znot}) when $z_i^\ell = 0$. 
Still, those constants should be preferably as small as possible.\footnote{For more detail,  see \nameref{app:relu_bounds} (in Appendix).}

As you may have guessed, or may have painfully experienced yourself, 
handling big--M constants is not very convenient: 
making them too large leads to numerical issues and making them too small leads to an invalid formulation. 
Many modern solvers support \emph{indicator constraints} for expressing 
big--M constraints and delegating the calculations involved to the solver. 
To formulate as in \cite{fischetti2018dnnmilp} but with the notation above, 
we replace constraints (\ref{eq:z}) and (\ref{eq:znot}) as follows: 
\begin{align}
z_i^\ell = 0 \rightarrow h_i^\ell \leq 0 \tag{\ref*{eq:z}'} \label{eq:z_ind}
\\
z_i^\ell = 1 \rightarrow \bar{h}_i^\ell \leq 0 \tag{\ref*{eq:znot}'} \label{eq:znot_ind}
\end{align}

\subparagraph{Modeling the Entire Network} 

Now we can effectively model the function $f(\vx)$ represented by a neural network using Mixed-Integer Linear Programming~(MILP) as follows:
\begin{subequations}
\begin{align}
\vx = \vh^0 \label{eq:x} \\
\begin{array}{r} (\ref{eq:W}) \\ ~ \end{array} & ~~ \begin{array}{ll} \forall & \ell \in \{1, \ldots, L+1\}, \\ & i \in \{1, \ldots, n_\ell\} \end{array} \label{eq:preact} \\
\begin{array}{r} (\ref{eq:g})\text{ to }(\ref{eq:bin}) \\ ~ \end{array} & ~~ \begin{array}{ll} \forall & \ell \in \{1, \ldots, L\}, \\ & i \in \{1, \ldots, n_\ell\} \end{array} \label{eq:act} \\
g^{L+1}_1 = y \label{eq:y}
\end{align}
\label{eq:nn_mip}
\end{subequations}
The whole-network formulation models the preactivation output $g_i^\ell$ of each neuron $i$ in layer $\ell$ in terms of its input $\vh^{\ell-1}$ with constraint (\ref{eq:preact}) 
and the output $h_i^\ell$ in terms of $g_i^\ell$ with constraint (\ref{eq:act}). 
Both constraints are applied to the $L$ hidden layers with ReLU activation, but only the first one is applied to the output neuron in layer $L+1$. 
We do not apply the ReLU activation to this last neuron because it would prevent us from representing negative outputs with $f(\vx)$.
The other two constraints, (\ref{eq:x}) and (\ref{eq:y}), are just ``syntax sugar'' between the internal naming of the network decision variables at each end, $\vh^0$ and $g_1^{L+1}$, and the external naming of decision variables in Figure~\ref{fig:nn_figure} and formulations (\ref{eq:f_cons}) and (\ref{eq:f_obj}), $\vx$ and $y$.  
Hence, 
we can model either formulations (\ref{eq:f_cons}) or (\ref{eq:f_obj}) 
as an MILP 
by denoting $f(\vx)$ as $y$ along with the constraints in~(\ref{eq:nn_mip}) for binding $\vx$ to $y$. 

\paragraph{What About Other Choices?} 

Readers who are more familiar with machine learning may have questions because they have seen more sophisticated architectures elsewhere, 
such as 
other activation functions that have gained  popularity recently, 
convolutional neural networks as a bare minimum in computer vision, 
and transformers as the standard architecture in contemporary natural language processing. 
Continuing with the longer answer, 
we will now discuss if and how other common architectures and variants are amenable for use in constraint learning as presented here. 

First, there are other activation functions that are also piecewise-linear, 
and as such they can also be modeled using a MILP formulation~\citep{vielma2010pwl}. 
Some of the most popular examples include 
the leaky ReLU~\citep{maas2013lrelu}, parametric ReLU~\citep{he2015prelu}, and maxout~\citep{goodfellow2013maxout}; 
with a formulation for the latter described in~\cite{serra2018bounding}. 

Second, 
there are some popular activation functions used to shape the network output as a probability, 
such as the sigmoid $\sigma(u) = \dfrac{1}{1 + e^{-u}}$~\citep{rumelhart1986learning} for single (scalar) outputs.\footnote{For probabilities over multiple (vector) outputs, see \nameref{app:softmax} (in Appendix).}  
Unlike the original setting that we described in Figure~\ref{fig:nn_figure} and refined in~\nameref{ap:nn_structure_map} (in Appendix), 
in this case we are also applying an activation function to the internal representation of the network output $g_1^{L+1}$, i.e., 
$y = \sigma^{L+1}\left(g_1^{L+1}\right) = \dfrac{1}{1 + e^{-g^{L+1}_1}}$. 
If the neural network models a constraint as in formulation (\ref{eq:f_cons}), 
i.e., $f(\vx) \leq d$ for $d$ constant, 
then we can reformulate the sigmoid output as a linear constraint: 
$f(\vx) = \dfrac{1}{1 + e^{-g_1^{L+1}}} \leq d$ implies that $g_1^{L+1} \geq \ln(d) - \ln(1-d)$. 
In turn, if the neural network models the entire objective function, as in formulation (\ref{eq:f_obj}) but with $\vc = \mathbf{0}$, i.e., the objective function is simply $\max f(\vx)$, 
then we can use the monotonicity of the sigmoid function to maximize (or minimize) directly on $g_1^{L+1}$. 
In other words,
we maximize $g_1^{L+1}$ instead of $y = \sigma^{L+1}(g_1^{L+1})$. 
Other possibilities, such as an objective of the form $\max \vc^T \vx + f(\vx)$ for $\vc^T \neq \mathbf{0}$, 
would require a more sophisticated approach, 
which is beyond the scope of this tutorial. 

Third, 
there are other types neural networks in which the same parameters are applied over multiple parts of the input. 
For example, that is the case of convolutional neural networks used in computer vision~\citep{lecun1998cnn}, 
where a filter consisting of a  small set of parameters is applied to multiple regions centered on different pixels of the image. 
This concept is generalized to other forms of structured data with graph convolutional networks~\citep{kipf2017gcn}, 
where a filter is applied to the neighborhood of each node of a graph representing the data. 
In both cases, 
every filter is equivalent to multiple neurons in terms of our setting in this tutorial, 
each of which applying the same transformation but to different regions of the input. 
Therefore, 
it is also possible to work in constraint learning with convolutional and graph convolutional neural networks with piecewise-linear activation functions. 
In fact, 
\cite{zhang2024gnn} 
extend the OMLT framework for that purpose. 
But because every filter corresponds to multiple neurons, 
this compact representation may 
disguise the much larger MILP formulations that are produced---to the point that it may become difficult to solve them, 
as we will see in \nameref{sec:how} 

Finally, and less fortuitously, 
some other activation functions that have gained popularity in recent years are not piecewise-linear, 
such as the Gaussian Error Linear Unit (GELU)~\citep{hendrycks2016gelu}, 
the Sigmoid-weighted Linear Unit (SiLU)~\citep{elfwing2018silu}, and 
Swish~\citep{ramachandran2018searching}. 
However, they are more often associated with generative models~\citep{shazeer2020glu} and with improving predictive performance in large-scale models.
Since one of the tacit premises of constraint learning is that of using the neural network as a more tractable approximation 
for nonlinear parts of a formulation, 
which entails the use of modestly-sized models and activation functions more amenable to optimization, 
those activations would be less appealing in this case.

\subsection{Where Have Neural Networks Been Used?}
\label{sec:where}

Neural networks have been used in many constraint learning applications, 
most notably to represent chemical engineering processes and power systems dynamics. 
We outline a non-exhaustive list of references in Table~\ref{tab:where} 
by approximate chronological order, briefly describing the type of optimization problem approached and the purpose of the neural network within that optimization model. 
References presenting multiple applications~\citep{schweidtmann2019global,dumouchelle2022neur2sp,dumouchelle2024neur2ro,dumouchelle2024neur2bilo,liu2025simulator} are separated with horizontal lines for ease of reading. 

\begin{table}[h!]
    \scriptsize
    \centering
    \begin{tabular}{@{}c@{}||@{}c@{}||@{}c@{}}
        \textbf{Reference} & \textbf{Type of problem(s)} & \textbf{Prediction(s) involved} \\ \hline \hline
        \cite{huang2003catalyst} & Catalyst composition & Catalytic performance \\
        \cite{fernandes2006fischertropsch} & Synthesis of a specific oil product & Outcome of a catalytic process  \\
        \cite{fahmi2012biodiesel} & Biodiesel production plant design & Unit operation, thermodynamics, mixing \\ 
        \cite{santanna2017psa} & PSA operation for N$_2$/CH$_4$  separation & Nitrogen/methane purity and recovery \\ 
        \cite{lombardi2017edml} & Computational workload dispatch & Processor efficiency for assignment \\
        \cite{say2017nonlinear} & HVAC, navigation, reservoir control & State transition in sensors \\ \hline
        \cite{schweidtmann2019global} & Operation of a fermentation process & Fermentation outcome \\ 
        & Operating point of a compressor plant & Compressor power consumption \\ 
        & Operating point of a cumene process & Process-unit outputs throughout the plant \\ \hline
        \cite{grimstad2019surrogate} & Offshore oil well operation & Well production, pipeline pressure drop \\ 
        \cite{bergman2022janos} & Scholarship allocation & Scholarship acceptance \\
        \cite{kody2022ac} & Unit commitment & AC power-flow solution \\
        \cite{murzakhanov2022dscopf} & AC power flow dispatch & Safisfaction of security constraints \\ \hline
        \cite{dumouchelle2022neur2sp} & Facility location & Customer demand \\
         & Investment decision & Return on investment \\
         & Server allocation & Client demand \\
         & Pooling & Product demand \\ \hline
        \cite{alcantara2023chance} & Energy production & Node voltage \\
        \cite{alcantara2023distributional} & Electricity pricing & Demand, day-ahead price, solar generation\\
        \cite{mcdonald2024molecular} & Molecular design & Molecule boiling point \\
        \cite{kronqvist2024approx} & Day-ahead power generation \& trade & Optimal intra-day generation \\ \hline
        \cite{dumouchelle2024neur2ro} & First-stage robust knapsack decision & Optimal second-stage objective \\
        & First-stage investment decision & Optimal second-stage objective \\ \hline
        \cite{dumouchelle2024neur2bilo} & Knapsack minimization for a leader & Optimal knapsack value for a follower \\
        & Network nodes to protect & Profit following network attack \\
        & Donation allocation & Return on projects developed by recipient \\
        & Network design & Total network travel time \\ \hline
        \cite{maragno2025miocl} & Diet planning & Diet palatability \\
        \cite{liu2025dca} & Data center allocation & Electricity price \\ \hline
        \cite{liu2025simulator} & Jet engine turbine blade design & Maximum blade temperature \\
        & Conversion of natural gas to H$_2$ & H$_2$ and N$_2$ concentrations, other process outputs \\ \hline
        \cite{feng2026windfarm} & Wind farm planning & Wind power production \\
        \cite{klimek2026aviation} & Kerosene production & Thermochemical and catalytic subprocesses \\
        \cite{luthje2026carnot} & Design and operation of Carnot batteries & Thermodynamic properties \\
        \cite{vahabzadeh2026hydrogen} & Underground hydrogen storage and operation & Hydrogen recovery factor \\
    \end{tabular}
    \caption{Applications from constraint learning studies in which neural networks are used.}
    \label{tab:where}
\end{table}

We have deliberately excluded from Table~\ref{tab:where} other optimization problems discussed in those and in other references if the optimization model embeds another type of predictive model; 
such cases are briefly outlined in \nameref{sec:not} 
For brevity, we have also excluded more recent references revisiting the same settings with a different approach. 

On the other hand, 
we have deliberately included some references in Table~\ref{tab:where} describing applications in which a neural network is used but not explicitly embedded in an optimization model; 
such cases correspond to the first four rows. 
All of them are applied optimization problems in chemical engineering, sharing the common theme of using a neural network as a surrogate for a complex system. 
Among them, the closest to the constraint learning setting is \cite{fahmi2012biodiesel}, 
which uses a disjunctive programming formulation.  
Most of them predate the popularization of the ReLU activation in neural networks, 
which was a turning point because the ReLU activation can be embedded in more tractable optimization models than the activation functions used earlier.

\section{When \emph{Not} to Use a Neural Network in Optimization Models?}\label{sec:not}

We have seen that  
a neural network can be a very expressive model for function approximation, 
which can be helpful in many applications. 
However, we may not need to always go that far. 
For example, there are many studies based on embedding much simpler predictive models into optimization models, such as the non-exhaustive list of references 
illustrated in Table~\ref{tab:not}. 
We list the types of models by  perceived complexity, 
citing earlier or most relevant references using each type.

\begin{table}[h!]
    \scriptsize
    \centering
    \begin{tabular}{@{}c@{}||@{}c@{}}
        \textbf{Type of predictive model} & \textbf{References} \\ \hline \hline
         Linear regression &  \cite{bertsimas2016chemo},  \cite{verwer2017auction} \\
         Logistic regression & 
         \cite{bergman2022janos} \\ 
         Support vector machine & \cite{hooshmand2025integration} \\
         Multivariate adaptive regression spline & \cite{martinez2017splines} \\
         Regression tree & \cite{verwer2017auction}, \cite{ammari2023trees} \\
         Decision tree & \cite{bonfietti2015embedding},  \cite{lombardi2017edml}, \cite{kudla2018oneclass} \\
         Random forest & \cite{bonfietti2015embedding} \\
    \end{tabular}
    \caption{Constraint learning studies in which other types of predictive models are embedded.}
    \label{tab:not}
\end{table}

Before training any neural network, 
it might be wise to run at least a linear regression and evaluate how good the predictions and decision quality can be with it. 
Moreover, as we will see next, the size of a neural network may affect how long it takes to solve the optimization model embedding it, 
sometimes in a prohibitive way if we are not careful. 
Hence, smaller and simpler can be better.

\section{How to Use Neural Networks in Optimization Models?}\label{sec:how} 

So far we have discussed the potential applications, expressive power, and modeling convenience of embedding neural networks in MILP models. 
But as we have just alluded to in the previous section, dimensioning the neural network properly is crucial for a successful use in constraint learning. 

\subsection{The Effect of Neural Network Design Choices}\label{sec:effect}

In order to choose more appropriate neural networks to embed in optimization models, 
we may benefit from insights from both the machine learning and the mathematical optimization communities. 
In what follows, 
we will discuss three dimensions that we believe are reasonably well understood by those communities, 
although what each community knows may not be fully realized by the other.

\paragraph{\textbf{Dimension 1:} Number of neurons.}\label{sec:number_neurons} 
The total number of neurons $N = \sum_{\ell=1}^L n_\ell$ 
is perhaps one of the best metrics to start discussing performance tradeoffs in constraint learning; but also one for which the priorities of machine learning and the mathematical optimization practitioners may not be fully aligned.

From a machine learning perspective, more neurons is better: 
it is intuitive that a neural network with more neurons is more expressive; and we can say even more if the activation function is the ReLU.  
Consider again the example in Figure~\ref{fig:pwa} and the explanation about linear regions in Page~\pageref{pag:linear_region}:  
we can think of a neural network with ReLU activations as a piecewise-linear function, 
where each piece corresponds to a different combination of neurons being active or inactive. 
If we have $N$ neurons, 
then we may potentially have a piecewise-linear function with $2^N$ pieces. 
Indeed, at some point we may also have too many neurons and the neural network may actually \emph{overfit} for its purpose, 
meaning that it performs well on the data used for training because of \emph{memorization} instead of \emph{generalization}, 
and by consequence it does not perform as well on unseen data. 
The neural network may also become too expensive to train or use at that point. 
However, as we will see next, 
we are unlikely to reach such a point in constraint learning 
because of other adverse effects. 

From a mathematical optimization perspective, 
less neurons is faster: 
it is intuitive that solving an optimization model embedding a larger neural network takes longer. 
If we consider the formulation~(\ref{eq:W})--(\ref{eq:bin}) mapping inputs to output of each neuron with ReLU activation, 
one of the main issues is the use of big M constraints (\ref{eq:z}) and (\ref{eq:znot}), 
for which the values of the constants $H^\ell_i$ and $\bar{H}^\ell_i$ must be as large as the variables $h^\ell_i$ and $\bar{h}^\ell_i$ could ever be for any valid input. 
At a microscopic level\footnote{For a broader discussion about alternatives to formulation~(\ref{eq:W})--(\ref{eq:bin}), see \nameref{ap:formulations} (in Appendix).}, 
the main issue with such big M constraints is what happens in the linear relaxation when those output variables take much smaller values. 
For example, 
if we relax the binary variable $z^\ell_i$ to be continuous between zero and one, 
then we can satisfy constraint (\ref{eq:z}) with $h^\ell_i = 1$ and $z^\ell_i = 0.001$ in a formulation where $H^\ell_i = 1000$; although $h^\ell_i > 0$ requires $z^\ell_i = 1$. 
When a solution with many binary variables having small values is obtained at a node of the branch-and-bound tree, 
we get very limited guidance on how to branch next. 
Because we may risk exploring as many as $2^N$ branch-and-bound nodes, 
solver runtime may induce a smaller choice of $N$ than training, use, or overfitting would. 

\paragraph{\textbf{Dimension 2:} Size of the input.}\label{sec:size_input} Although less evident, 
the size of the input $n_0$ of a neural network is an important metric; 
and one for which 
a smaller size can be advantageous from both perspectives. 
Even if the size of the input is not typically under our control, 
understanding its effect helps explaining why 
neural networks that are otherwise identical yield very different results. 
Intuitively, a dataset with fewer dimensions may require a smaller neural network, 
but that is just the tip of the iceberg. 

From a machine learning perspective, fewer inputs means a simpler model: 
the neural network is modeling a function on a smaller domain. 
When we use the ReLU activation function, 
that surprisingly translates to modeling a piecewise-linear function with fewer pieces. 
Consider one more time the example in Figure~\ref{fig:pwa}, 
where a neural network with 5 neurons models a function with 12 pieces.  
If every piece, also known as a linear region, 
is defined by a combination of which neurons are active or not, 
then we could potentially have as many as $2^5=32$ regions. 
However, for a neural network with two inputs, 
we can only define up to 16 regions with five neurons in a single hidden layer; no matter how we change the network parameters\footnote{For a geometric argument, see \nameref{ap:geometry} (in Appendix).}. 
Similarly, a small network input also affects the number of regions defined by neurons after the first layer and, more generally, 
any hidden layer being too narrow also affects the number of regions defined by the layers that follow it~\citep{huchette2026polyhedral}. 
\cite{serra2018bounding} coined it as the \emph{bottleneck effect}.

From a mathematical optimization perspective, 
fewer linear regions means a smaller feasible set: 
if every linear region corresponds to a different assignment to the binary variables in $\vz$, 
then having fewer linear regions means that more of the binary assignments to $\vz$ are infeasible. 
When it comes to solving the MILP formulation by using branch-and-bound, 
we are likely to reach infeasible nodes more often with a neural network having fewer inputs in the formulation. 
Hence,
even if the formulation is weak because of the big M constraints, 
we may not need to go very deep before coming across dead ends with the branch-and-bound algorithm. 
For the same reason, 
local search methods that explore their neighborhood more thoroughly---and thus resemble more the enumerative approach of the branch-and-bound algorithm---tend 
to be more competitive when the neural network has smaller input size. 
For smaller inputs, 
the method by \cite{perakis2022sampling} tends to produce better solutions in comparison to the method by \cite{tong2024lp}, 
and the latter tends to produce better solutions in comparison to the method by \cite{tong2026gradient}. 
Conversely, 
we can find better solutions within a fixed time limit for optimization models embedding larger neural networks, and especially neural networks with larger input size, 
if we use heuristics that move faster from solution to solution; 
hence reversing the order of preference among the heuristics cited above. 
Moreover, 
we encourage the reader to be less afraid of trying larger neural networks when the input is smaller simply because the solver runtime is not going to be impacted as much.

Yet from a modeling perspective, 
this means that identifying any form of dimensionality reduction can be beneficial if we would like to use an MILP solver. 
Part of the modeling process consists of identifying what would be the inputs of the neural network, 
each input also being known as a \emph{feature} among machine learning practitioners. 
For example, 
if we observe that two features are highly correlated, 
meaning that one can be a proxy for the other, 
then there might be also a benefit to solver runtime if one of those were to be dropped from the neural network that will be trained.

\paragraph{\textbf{Dimension 3:} Connectivity between neurons.}\label{sec:connect} 
We often observe neural networks being measured and compared in terms of their number of \emph{parameters}. 
For example, the number of parameters is often used as a proxy for the capability of generative models: 
when a new and arguably more capable model is announced, its larger number of parameters is often emphasized. 
Indeed, 
more capable neural networks tend to have larger dimensions, and therefore it seems intuitive that having a larger number of parameters should be better. 
However, 
reducing the number of parameters of a neural network can be beneficial. 
This is usually done with \emph{network pruning}, 
where we fix some of the parameters, typically the weights, to zero after training and keep them at zero going forward.

From a machine learning perspective, 
the training of larger neural networks is often more successful for at least a couple of reasons.  
First, training is solved as an unconstrained optimization problem using first-order methods, 
where the objective function is known as the \emph{loss function}; 
and larger neural networks define a smoother loss, yielding better convergence~\citep{li2018landscape,ruoyu2020landscape}. 
Second, having more neurons per layer prevents the layers from becoming entirely inactive and stall training~\citep{riera2022jumpstart}. 
But after training, 
we can remove as much as half of the parameters of sufficiently large neural networks and still recover or improve upon the original predictive performance by retraining after pruning, which is known as \emph{finetuning}~\citep{hoefler2021sparsity}; 
where ``sufficiently large'' depends on what we train the neural network for~\citep{liebenwein2021lost}, 
and finetuning is generally brief because the loss landscape from zeroing out pretrained networks yields a faster convergence~\citep{jin2022generalization}. 
But curiously, and also fortunately, the number of linear regions is not affected by a moderate amount of pruning~\citep{cai2023pruning}. 
Finally, a simple heuristic known as \emph{magnitude pruning}, which consists of removing the weights with smallest absolute value, tends to be very effective for moderate pruning~\citep{pham2026sparse}. 

From a mathematical optimization perspective, 
MILP solvers have been historically designed to exploit the sparsity in models. 
\cite{dey2018sparse} observe that the number of nonzero coefficients in constraints is often very small in benchmark problems, and that 
Linear Programming~(LP) solvers rely on linear algebra routines that can be made more efficient by exploiting matrix sparsity. 
One such example is the study by \cite{crowder1983large} that first solved several large-scale but sparse LP models. 
Validating how this trend continued to influence MILP solvers design,  
\cite{matthias2014sparsity} empirically verified a significant slowdown when models are augmented with a few dense constraints. 
In fact, 
some of the applications previously mentioned in \nameref{sec:where} have deliberately applied network pruning to exploit the benefit of sparsity in MILP solvers~\citep{say2017nonlinear,kody2022ac} 
and recent methodological studies have considered this topic in more detail~\citep{cacciola2024structured,pham2026sparse,lastrucci2026pruning}.

\subsection{A Case Study}

Having discussed mathematical formulations, potential applications, and hypothetical design choices of constraint learning, 
here we stand, poor fools, 
barely wiser in practice than before!\footnote{This is a pastiche of a quote in \cite{faust} that translates as 
\emph{``I HAVE, alas! Philosophy,
Medicine, Jurisprudence too,
And to my cost Theology,
With ardent labour, studied through.
And here I stand, with all my lore,
Poor fool, no wiser than before.''}} 

We address this gap with a computational experiment. 
We use a hypothetical problem with a synthetic dataset generated from a function of our choosing 
because that allows us to compare neural network predictions with actual outcomes (i.e., the \emph{ground truth}). 
This is typically not necessary when neural networks are used for prediction because we can reserve part of the dataset (the \emph{validation set}) for comparing different neural networks. 
However, the solution coming out of the optimization model may assign inputs to the neural network that are outside the dataset.

In this example, 
we will focus on two factors mentioned in \nameref{sec:effect}, (i) \emph{number of neurons} and (ii) \emph{connectivity between neurons}, 
illustrating how they affect prediction quality, decision quality, and computational runtime. 
Those factors often require empirical calibration in constraint learning. 
The other factor, \emph{size of the input}, 
would be harder to illustrate with a single example without an extensive discussion about dimensionality reduction.  

Some technical details are intentionally omitted from the description below for ease of reading. 
For more details, we refer the reader to \nameref{app:dataset_generation} and to \nameref{ap:experiments} (both in Appendix). 
Moreover, we note that this experiment is simply intended to substantiate our discussion, 
hence why we are focusing on a single problem.

\paragraph{The Problem}
A firm must decide how to spend their marketing budget for advertising a product across eight channels: television; radio; social media, including which platform to use; influencer sponsorship; billboard advertising; browser advertising; email advertising; and event sponsorship. 
Those decisions are influenced by product price, sales in the last month, and the marketing budget of the largest competitor; 
which we assume as \emph{exogenous parameters} (i.e., not under our control). 
The goal is to maximize the estimated revenue from advertising decisions based on historical data associating past advertising decisions and exogenous parameters with product sales.

\paragraph{The Optimization Model} 
To describe a model for this problem, we first introduce some notation.
Problem dimensions: 
\begin{itemize}
    \item $n_c = 8$: the number of advertising channels.
    \item $n_p = 5$: the number of social media platforms.
    \item $n_e = 3$: the number of exogenous parameters.
\end{itemize}
Bounds on the problem data:
\begin{itemize}
    \item $l^c_i \in \mathbb{R}_+$: lower bound on spending in channel $i$ for $i = 1, \ldots, n_c$.
    \item $u^c_i \in \mathbb{R}_+$: upper bound on spending in channel $i$ for $i = 1, \ldots, n_c$.
    \item $l^e_i \in \mathbb{R}_+$: lower bound of exogenous parameter $i$ for $i = 1, \ldots, n_e$.
    \item $u^e_i \in \mathbb{R}_+$: upper bound of exogenous parameter $i$ for $i = 1, \ldots, n_e$.
\end{itemize}
Decision variables:
\begin{itemize}
    \item $x_i$: how much the firm spends in channel $i$ for $i=1, \ldots, n_c$.
    \item $z_i$: if social media platform $i$ is used for $i = 1, \ldots, n_p$.
\end{itemize}
Instance data:
\begin{itemize}
    \item $B \in \mathbb{R}_+$: marketing budget for advertising the product.
    \item $\xi_i \in \mathbb{R}_+$: observed value for exogenous parameter $i$ for $i = 1, \ldots, n_e$.
\end{itemize}
Now we can formulate this problem as follows:
\begin{subequations}\label{mip:case}
\begin{align}
\max ~~ & f(\vx, \vz, \boldsymbol{\xi}) \\
\text{s.t.} ~~ & \sum_{i=1}^{n_c} x_i \leq B \label{eq:case_begin} \\
& \sum_{i=1}^{n_p} z_i = 1 \\
& \vl^c \leq \vx \leq \vu^c \\
& \vz \in \{0, 1\}^{n_p} \label{eq:case_end}
\end{align}
\end{subequations}
Where function $f : \mathbb{R}_+^{n_c + n_p + n_e} \rightarrow \mathbb{R}$ models revenue based on the values assigned to the decision variables in vectors $\vx$ and $\vz$ as well as on the values of the exogenous parameters in vector $\boldsymbol{\xi}$.

\paragraph{The Dataset} 
As a reminder, 
in constraint learning 
we work under the assumption that we cannot optimize directly over $f$ even if we know the exact expression for $f$, 
which is why we replace it by a \emph{proxy function} $\hat{f}$ learned with a neural network. 
Hence, 
we constructed a nonlinear function $f$ intended to capture diminishing returns from spending more in each channel, different revenue peaks based on what social media platform is used, and declining revenue with price and competitor spending. 
Then we generated a synthetic dataset with observations of the form $\big(\bar{x_1}, \bar{x_2}, \ldots, \bar{x_{n_c}}, \bar{z_1}, \ldots, \bar{z_{n_p}}, \bar{\xi_1}, \ldots, \bar{\xi_{n_e}}, f(\bar{\vx},\bar{\vz},\bar{\boldsymbol{\xi}}) + \varepsilon\big)$, 
where $\varepsilon$ is a small fluctuation based on random noise. 
This dataset is partitioned into training and validation sets, 
which are used for obtaining and evaluating neural networks to be used as a proxy function $\hat{f}$ in the place of $f$. 
The validation set is only used to evaluate neural network predictions. 
For evaluating downstream decisions based on a neural network acting as a proxy, we use the ground truth function $f$ directly. 
Other details can be found in~\nameref{app:dataset_generation} (in Appendix).

\paragraph{Neural Networks Evaluated}
Using 10 different random seeds to initialize the network parameters,  
we trained fully connected neural networks with $L=2$ hidden layers, 
each layer $\ell \in \{1, 2\}$ having width $n_\ell = w$. 
We trained small neural networks with $w = 10, 20, \ldots, 100$  neurons 
as well as large neural networks\footnote{We regard neural networks with width $w=100$ as being both small and large. Using $w=100$ as a baseline by including those results in both groups makes it easier to compare across groups. 
But if we need to choose one group only, we regard them as small.} with $w = 100, 200, \ldots, 1,000$ neurons. 
To each trained network, 
we applied magnitude pruning to remove a fraction $p$ of the weights, 
with \emph{pruning rate} $p = 0.0, 0.1, \ldots, 0.9$ 
and, except for $p=0.0$, 
finetuned the pruned network. 
Given the number of seeds, widths, and pruning rates used in the computational experiment, we trained a total of 1,900 neural networks. 
We use the number of weights as one of the metrics for comparing these neural networks in terms of size. The number of weights is calculated as 
\begin{equation}
    W(d,w,p) = (d w + w^2 + w) (1 - p)
\end{equation}
for an input size $d$, width $w$ in the two hidden layers, and pruning rate $p$.
Other details can be found in~\nameref{ap:experiments} (in Appendix).

\begin{figure}[b!]
    \centering
    \includegraphics[width=0.45\linewidth]{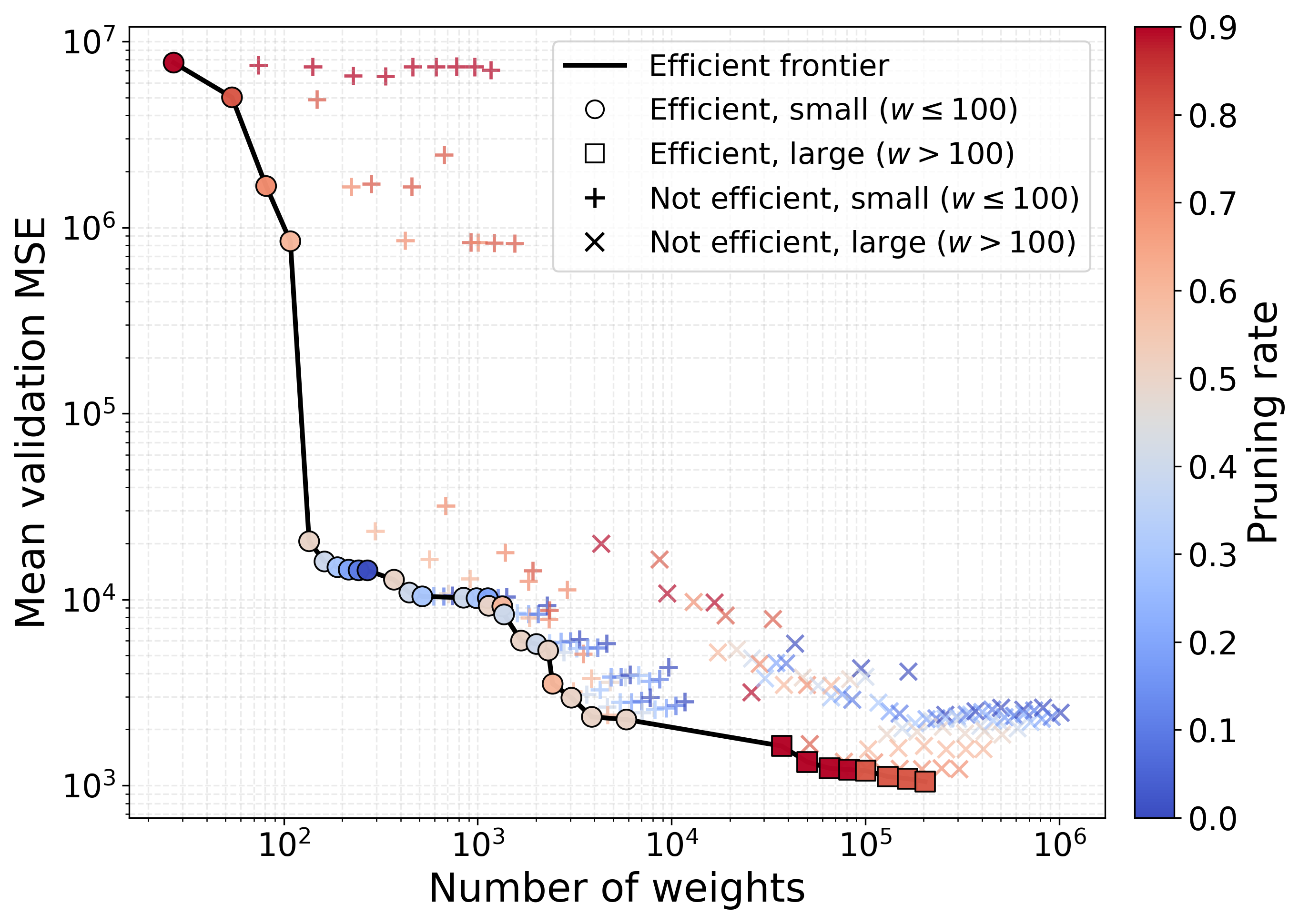}
    \caption{Comparison of predictive quality with number of weights of configurations (as a scatter plot), contrasted by pruning rate (as a color gradient). The configurations that are not dominated by any other in both predictive quality and number of weights define the efficient frontier (as a line chart).}
    \label{fig:pareto_frontier}
\end{figure}

\paragraph{Predictive Accuracy} 
Does the number of weights tell us \emph{everything}, or at least \emph{something}, about predictive accuracy?\footnote{We use words such as \emph{accurate} and \emph{accuracy} (as adjectives) to discuss the quality of predictions of different neural networks for regression, not classification. Hence, we are \emph{not} referring to the metric \emph{accuracy} (as a substantive) that is used in classification.}
We start by comparing (a) the number of weights $W(d,w,p)$ for different choices of width $w$ and pruning rate $p$ with (b) predictive performance in Mean Squared Error (MSE) on the validation set, averaging across random seeds. 
Indeed, we observe a size--accuracy tradeoff: 
the configurations that are not dominated by any other in both criteria define an \emph{efficient frontier},  
within which no other configuration is more accurate if smaller or less accurate if larger. 
Figure~\ref{fig:pareto_frontier} illustrates all 190 $(w,p)$ \emph{configurations} evaluated, 
where the efficient frontier contains 34 of those configurations.   
Although the most accurate configurations have a larger number of weights (see right side of the efficient frontier), 
note how those configurations at the efficient frontier tend to have higher pruning rates (see the markers becoming more red than blue). 
In other words: 
even if predictive performance improves with the number of weights, 
the most accurate models are pruned.  
Hence, there is a connection between the two metrics,
but this connection is insufficient as a proxy for network accuracy: 
the configurations with the largest number of weights, 
which have almost 10 times more weights than any configuration at the efficient frontier, 
perform in line with small networks.  
So we turn to network pruning to understand the nuance of this relationship.

\begin{figure}[b]
    \centering
    \includegraphics[width=0.45\linewidth]{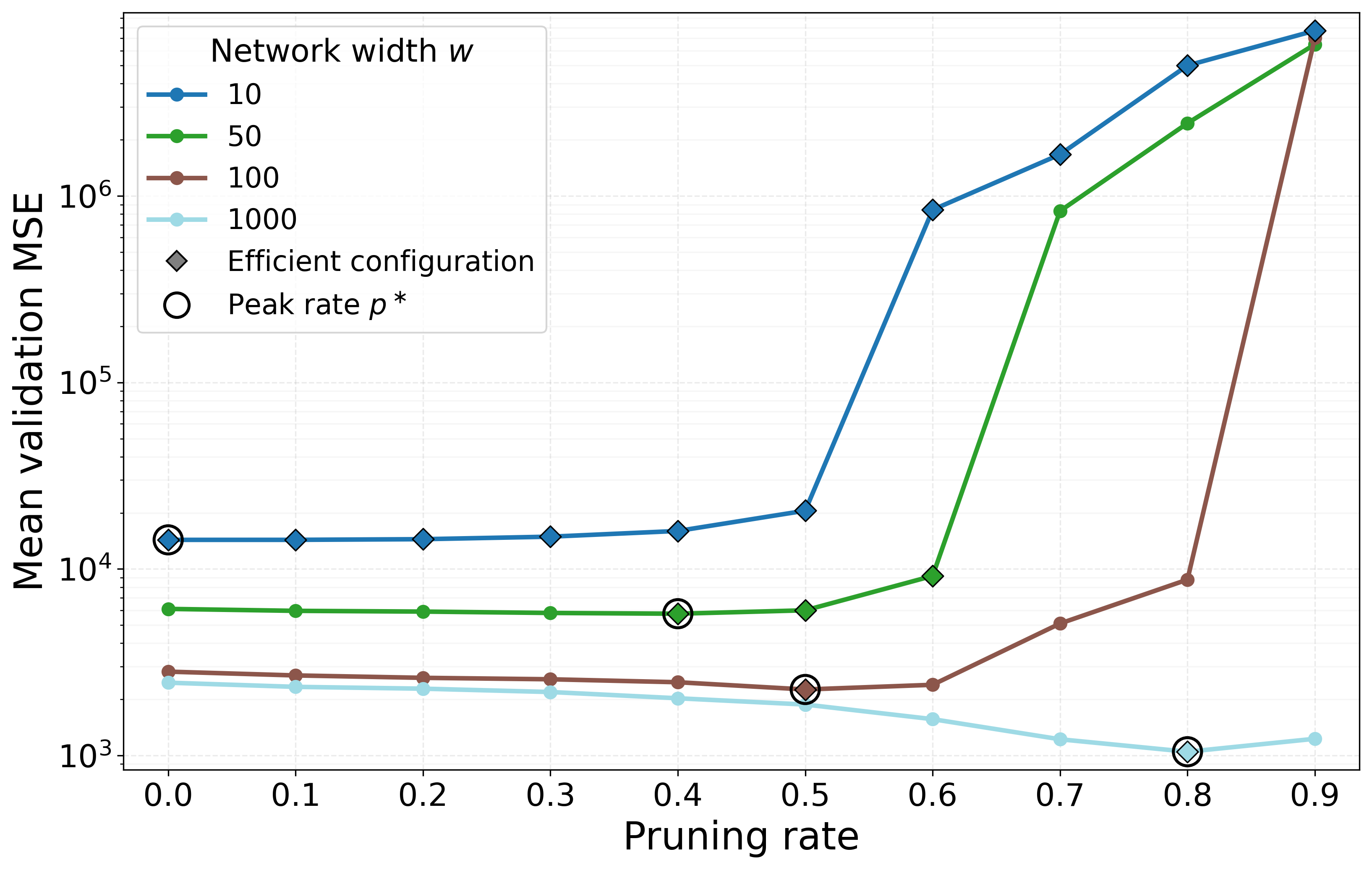}
    \caption{Comparison of predictive quality with pruning rate $p$ (as a line chart for each width $w$) for neural networks with width 10, 50, 100, and 1,000. The pruning rate $p*$ with smallest validation error for each width $w$ is circled. The configurations at the efficient frontier are diamond-shaped.}
    \label{fig:pruning_vs_mse_selected}
\end{figure}

\paragraph{Effect of Pruning on Accuracy} 
Can we always predict the effect of pruning in accuracy, 
or at least understand how it changes as neural networks get larger?\footnote{Cue: whenever a sentence starts with a question like this, the answer leans towards the what is expressed at the end of the question.}
We refine our understanding about this dichotomy between (i) more accurate networks being larger and (ii) the most accurate networks being pruned 
by directly comparing, 
for networks with each width $w$, the pruning rate $p$ with their predictive performance. 
Figure~\ref{fig:pruning_vs_mse_selected} presents line chars for networks with width 10, 50, 100, and 1,000, 
with each configuration $(w,p)$ at the efficient frontier being circled. 
By presenting the results in this way, 
we can take advantage of the following: 
if a line chart contains points at the efficient frontier, 
then the pruning rate $p^*$ with most accurate networks is the leftmost point. 
Figure~\ref{fig:pruning_vs_mse} (in Appendix) presents the results across all widths, 
to which we may refer for generality. 
From the plots, we note that $p^*$ tends to increase with $w$. 
In fact, 
Figures~\ref{fig:pruning_vs_mse_selected} and \ref{fig:pruning_vs_mse}  
show three pattens regarding which pruning rates are at the efficient frontier, depending on $w$, 
and those patterns help us understand how the pruning rate $p$ vary along the efficient frontier in Figure~\ref{fig:pareto_frontier}:
\begin{itemize}
    \item When the neural network dimensions are very small, any reduction in the number of weights worsens network accuracy, and so pruning is not beneficial. 
    That happens in Figure~\ref{fig:pruning_vs_mse} for $w=10$, 
    where pruning more implies lower accuracy, 
    and therefore all choices of $p$ are on the efficient frontier. 
    In Figure~\ref{fig:pareto_frontier}, that explains 
    the leftmost points (up to $3 \times 10^3$ weights), going from red to blue ($p=0.9$ to $p=0.0$) with the number of weights.
    \item For other small networks, the neural network is slightly more accurate with a low pruning rate $p^*$, and so moderate pruning is beneficial. 
    That happens in Figure~\ref{fig:pruning_vs_mse} for $w=20$ to $w=100$, 
    where some choices between $p=0.2$ and $p=0.6$  are on the efficient frontier. 
    In Figure~\ref{fig:pareto_frontier}, that explains 
    the middle segment of the efficient frontier (from $3 \times 10^3$ to $10^4$ weights) 
    going from red to neutral ($p=0.0$ to $p=0.5$), albeit with some fluctuations (such as $p=0.6$ in between).
    \item For large networks, the neural network becomes pronouncedly more accurate with a high pruning rate $p^*$, and so aggressive pruning is beneficial. 
    That happens in Figure~\ref{fig:pruning_vs_mse} for $w=500$ to $w=1,000$, 
    where only choices with $p=0.8$ and $p=0.9$ are on the efficient frontier. 
    Although the most beneficial pruning rate $p^*$ keeps increasing, 
    it never reaches the largest choice ($p^* \leq 0.8$ for $w \leq 1000$). 
    In Figure~\ref{fig:pareto_frontier}, that explains 
    the final segment of the efficient frontier (from $10^4$ weights onward) 
    having a drastic shift towards red ($p=0.9$) 
    while the largest \emph{parameter-efficient}\footnote{Although we only prune weights and thus we measure and compare number of weights across configurations, we use the term parameter-efficient for referring to the networks with best prediction for their size because that term is more popular.} network configuration (the rightmost point) has a slightly lower, but still high, pruning rate ($p=0.8$). 
\end{itemize}
From a predictive accuracy perspective, 
network pruning is therefore not only beneficial, 
but actually necessary to obtain large neural networks that are more parameter-efficient. 
However, we have not yet put that in the context of solving optimization models in  constraint learning.

\paragraph{Decision Quality and Optimization Bias} 
If neural networks with more weights and adequate pruning are better for prediction, 
are they also better for optimization? 
To answer this question, 
we carry out the following sequence of steps with each neural network that we trained or pruned: 
\begin{enumerate}
    \item We solve the problem with the neural network by maximizing $\hat{f}(\vx, \vz, \boldsymbol{\xi})$ subject to constraints (\ref{eq:case_begin}) to (\ref{eq:case_end})\footnote{We fix the values of the exogenous parameters in vector $\boldsymbol{\xi}$ and of the budget $B$ in constraint (\ref{eq:case_begin}) before solving the optimization model associated with this step. For more information, see \nameref{ap:experiments} (in Appendix).}. 
    That gives us a solution $(\vx^*,\vz^*)$ with the \emph{network optimal} value $\hat{f}(\vx^*, \vz^*, \boldsymbol{\xi})$. 
    \item We measure how good that solution actually is through the ground truth value $f(\vx^*, \vz^*, \boldsymbol{\xi})$.
    \item We measure the disagreement between the two by simply taking their difference, which \cite{lombardi2017edml} defines as the \emph{optimization bias} $\beta := \hat{f}(\vx^*, \vz^*, \boldsymbol{\xi}) - f(\vx^*,\vz^*, \boldsymbol{\xi})$. 
\end{enumerate}
When solving the optimization models, we set a time limit of 1,000 seconds to avoid spending too long with configurations that may not be the most parameter-efficient. As we will see from the plots, that time limit was sufficient for our experiments.

\begin{figure}[h]
    \centering
    \includegraphics[width=0.55\linewidth]{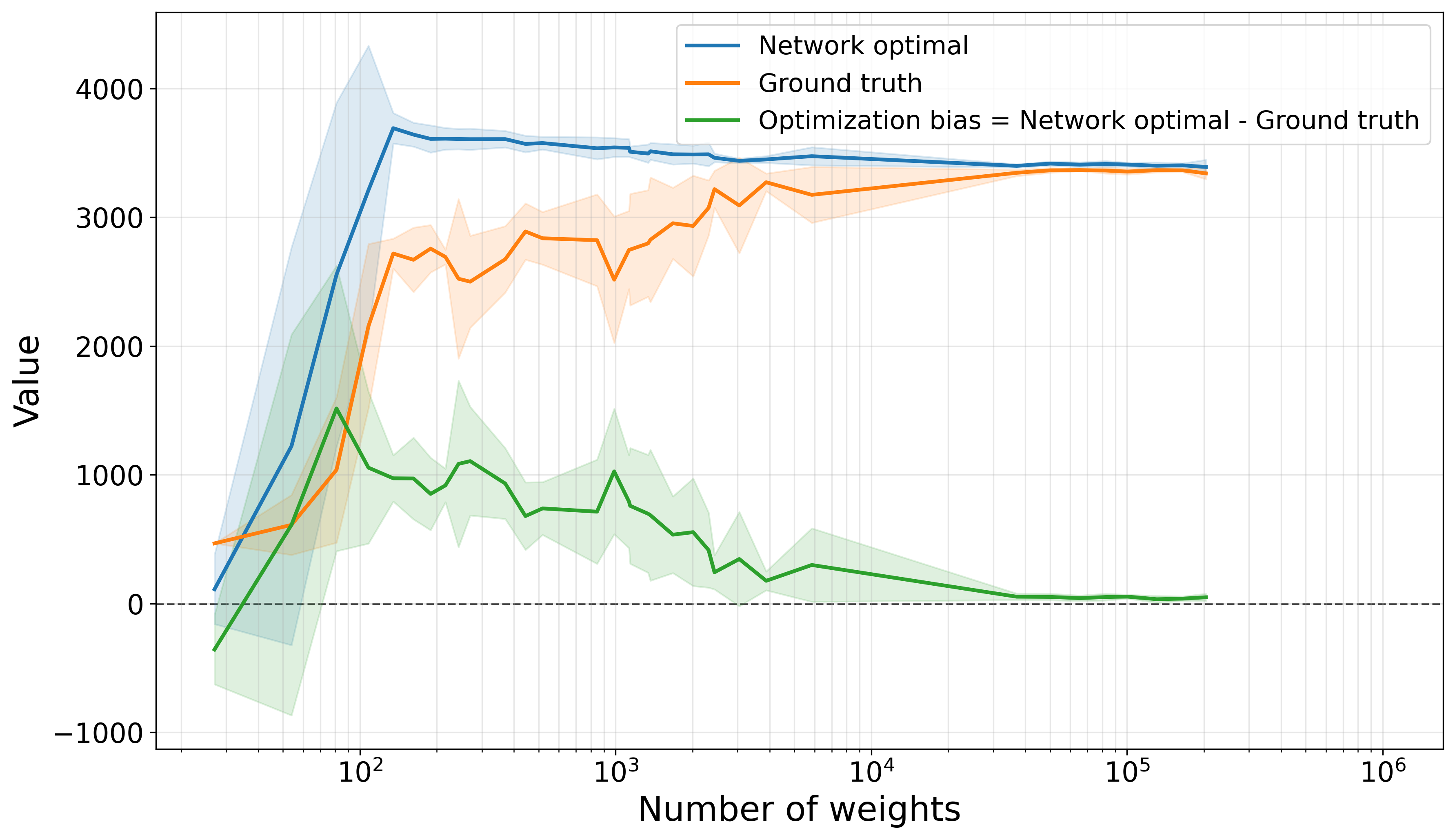}
    \caption{Comparison of decision quality with number of weights for configurations at the efficient frontier (as line charts with error bars). Using the solution obtained with the optimization model, we evaluate it as an input of the embedded network and compare it with the ground truth.}
    \label{fig:downstream_quality}
\end{figure}

Based on the values that we observe for configurations at the efficient frontier in Figure~\ref{fig:downstream_quality}, we can draw some conclusions for sufficiently large neural networks.  
First, and not very surprisingly, 
larger networks  naturally yield better solutions because the ground truth value tends to increase and the optimization bias tends to decrease. 
Second, and as expected from \cite{lombardi2017edml}, 
the optimization bias is typically positive in this case because 
the optimal solution from the neural network tends to be a solution that the neural network overestimates in maximization problems. 
Finally, and a little more subtly, 
the network optimal value can be useful on its own to evaluate decision quality because  
it slowly decreases and approaches the ground truth value in larger networks. 
Since the network optimal is the only of those three values that we would always have in practice,\footnote{When the ground truth is unknown, we cannot calculate the optimization bias, and thus we can only rely on the network optimal.} 
we may use it for an educated guess about the benefit of using even larger networks.

\begin{figure}[t]
    \centering
    \includegraphics[width=0.75\linewidth]{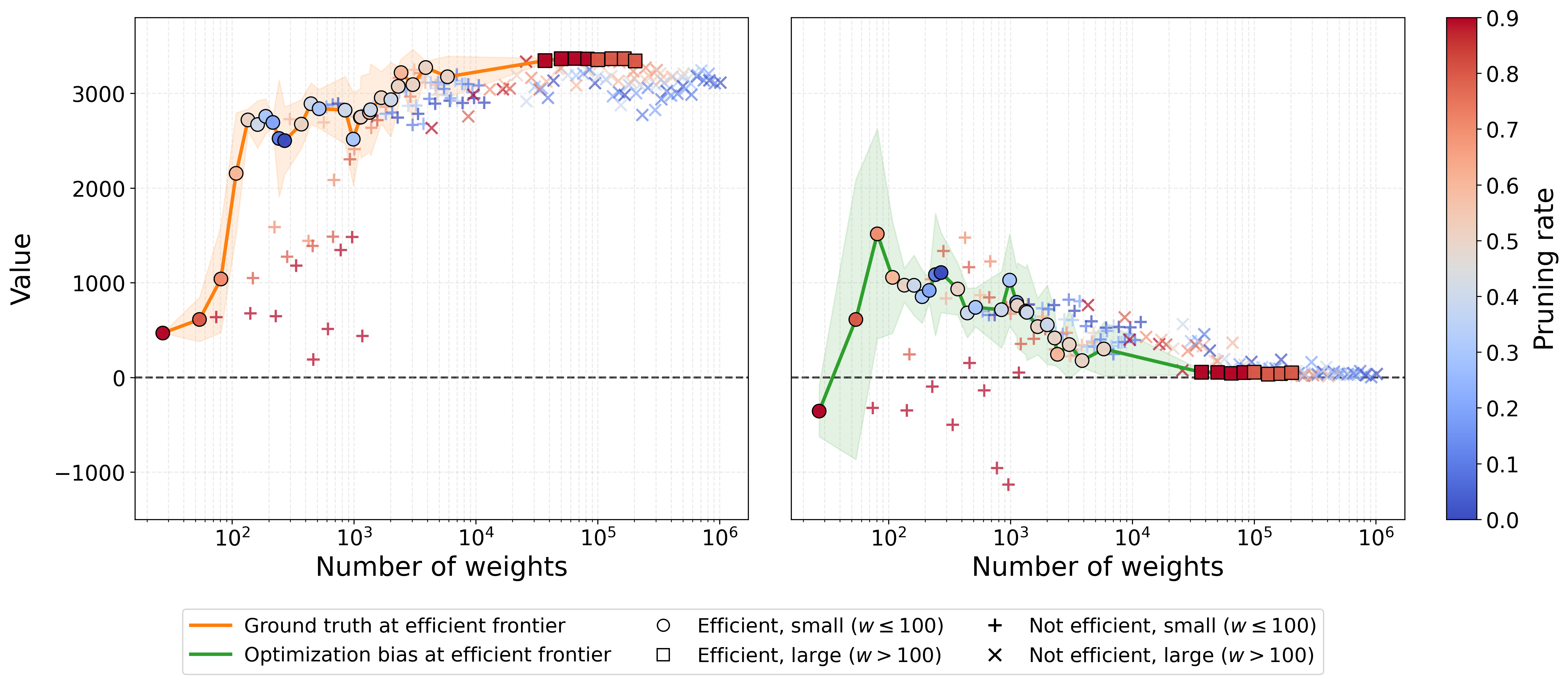}
    \caption{Comparison of decision quality with number of weights for configurations at the efficient frontier (as a line chart; the same ones from Figure~\ref{fig:downstream_quality}) vs. other configurations (as a scatter plot), contrasted by pruning rate (as a color gradient).
    For ground truth value, 
    104 configurations are below the curve defined by the efficient frontier, 15 are above, and 37 are outside the range.
    For optimization bias, 
    84 of those configurations are above the curve, 35 are below, and 37 are outside the range.}
    \label{fig:frontier_nonfrontier_quality}
\end{figure}

\begin{figure}[h]
    \centering
    \includegraphics[width=0.45\linewidth]{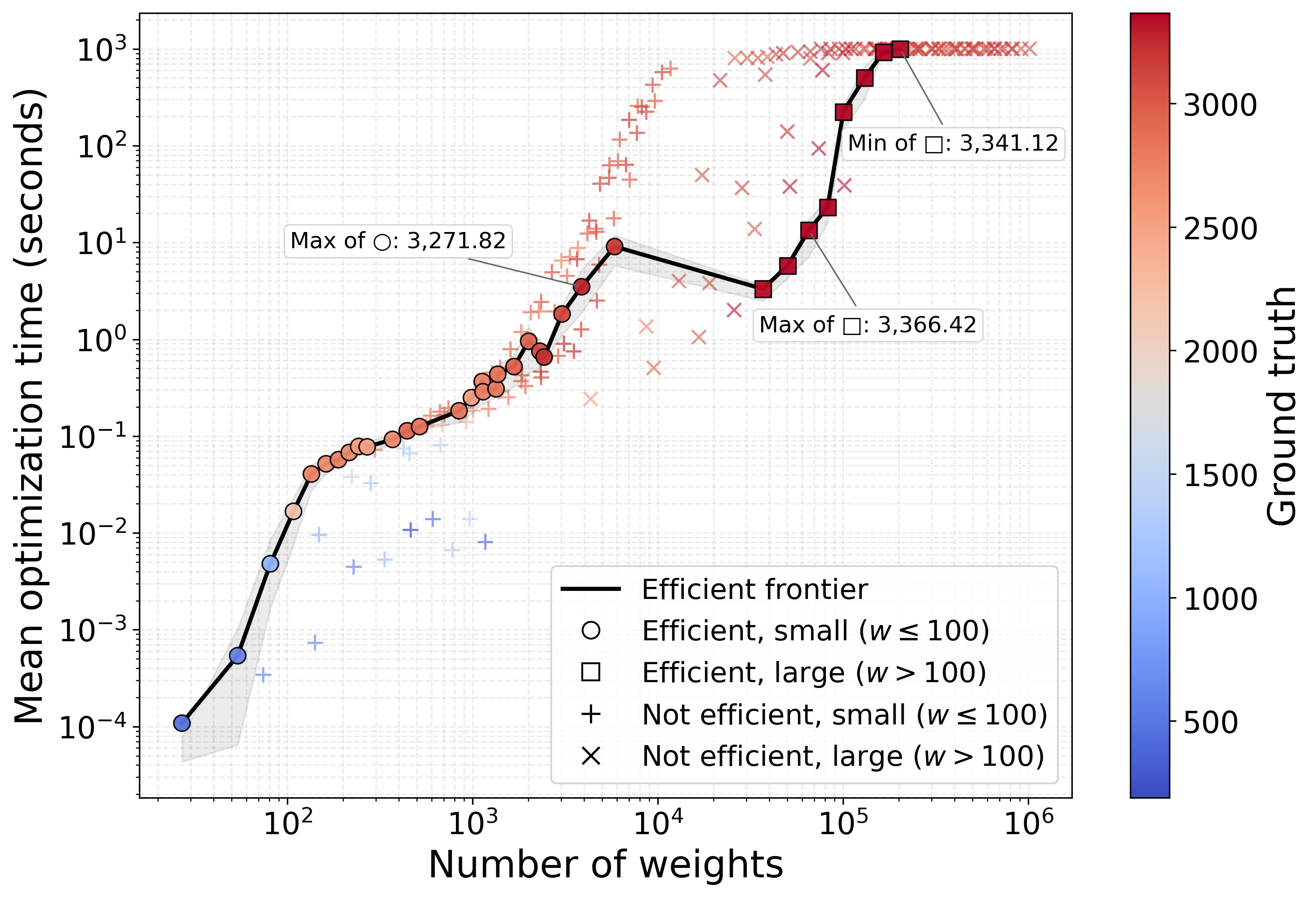} 
    \caption{Comparison of solution time with number of weights at the efficient frontier (as a line chart) vs. other configurations (as a scatter plot), contrasted by ground truth values (as a color gradient).}
    \label{fig:runtime_vs_size}
\end{figure}

Moreover, 
Figure~\ref{fig:frontier_nonfrontier_quality} shows that configurations at the efficient frontier are also more effective for optimization: 
the ground truth values for other configurations are often smaller or at best have comparable values if better, despite the optimization bias oscillating above and below the curve defined by the efficient frontier.

\paragraph{Comparing Runtime and Decision Quality} 
While the benefit in terms of both accuracy and decision quality diminishes when neural networks get sufficiently large, 
what happens to solver runtime? 
Figure~\ref{fig:runtime_vs_size} shows a steady 
increase in solver runtimes as networks get large. 
The diminishing benefit of using more solver runtime with a larger network becomes more pronounced 
if we consider the small difference in ground truth value across the large networks at the efficient frontier.

\begin{figure}[h]
    \centering
    \includegraphics[width=\linewidth]{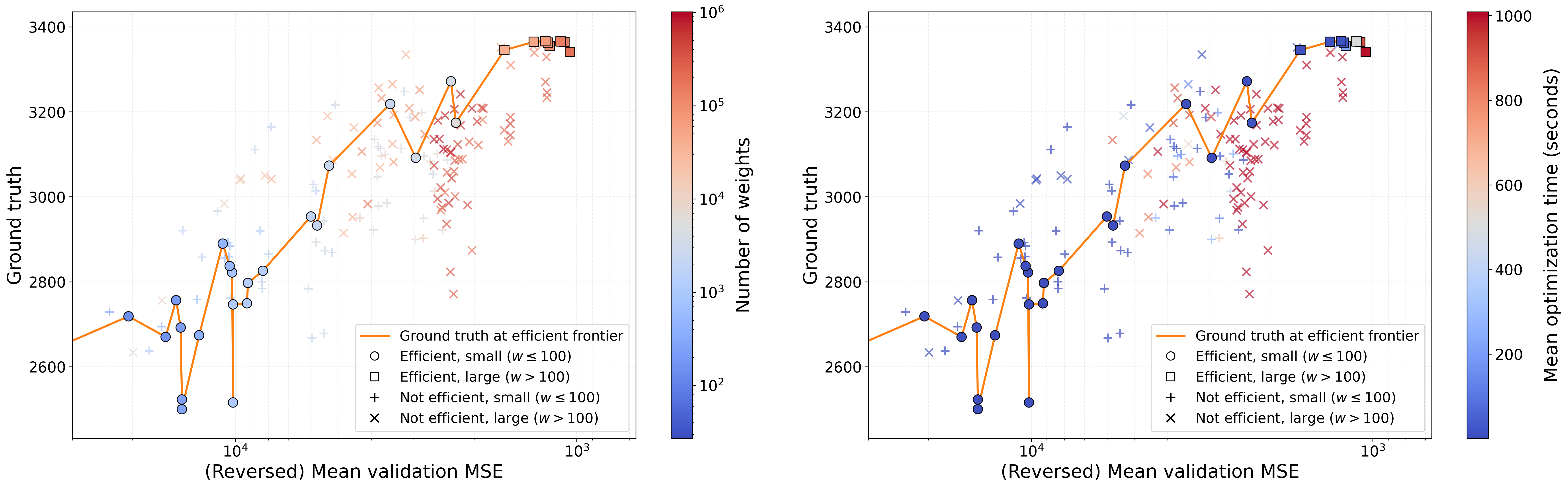}
    \caption{Comparison of ground truth value with predictive accuracy at the efficient frontier (as a line chart) vs. other configurations (as a scatter plot), contrasted on the left by the number of weights and on the right by solution time (as a color gradient in both; the plots are otherwise identical). The plots exclude configurations with mean validation MSE above $3 \times 10^4$. The correlation across for the points shown at the efficient frontier is -0.9535, dropping to -0.7885 for all points. }
    \label{fig:prediction_decision}
\end{figure}

\paragraph{Predictive Accuracy vs. Decision Quality} 
For a complementary view, what would happen if we primarily used accuracy for selecting which neural networks to embed? 
Figure~\ref{fig:prediction_decision} shows that there is indeed a positive correlation between predictive accuracy and decision quality, meaning a negative correlation between the log of mean validation MSE (reversed in the plot) and ground truth value, 
but also that this correlation does not tell the full story. 
By repeating this plot twice with the color gradient contrasting by different metrics, 
we can also dissuade the reader about why the configurations corresponding to points above the line chart defined by the efficient frontier are not advisable: 
those are networks with more weights and with longer solution time 
for either the same validation (i.e., by fixing the $x$ axis) or for the same ground truth value (i.e., by fixing the $y$ axis).

\section{Other Topics in Constraint Learning and at the AI--Optimization Interface}\label{sec:etc}

So far in this tutorial, we have only covered the fundamentals of constraint learning with neural networks. 
Inspired by a Sesame Street book, 
we have intentionally left this ``Everything Else''\footnote{At the end of Grover's visit to the Everything in the Whole Wide World Museum~\citep{stiles1974grover}, 
there is a door labeled ``Everything Else'' that leads to the outside. That is where Grover finds all the things that he did not see inside the museum.} section at the end to (i) provide references for 
 additional topics; and (ii) briefly discuss other prominent forms of integrating optimization with machine learning and Artificial Intelligence~(AI). 

We begin with directions that have been explored to improve or expand constraint learning:
\begin{enumerate}[(i)]
    \item We can reduce the solver runtime for the MILP formulation described in the tutorial 
if we can use smaller values for its big M constants, often called \emph{activation bounds}.  
Because finding the smallest possible value becomes as expensive as solving the entire model, 
there is extensive work on developing cost-effective bounding methods. 
Some of this work is described in~\nameref{app:relu_bounds} (in Appendix).
    \item We can use entirely different MILP formulations for embedding neural networks, 
    or at least generate cutting planes based on the same mechanism yielding such formulations. 
    Those formulations are often much larger but are also tighter, meaning that they may require exploring fewer branch-and-bound nodes before concluding that an optimal solution has been found. 
    Some of those formulations are mentioned in~\nameref{ap:formulations} (in Appendix).
    \item We can restrict the neural network to learn only convex proxy functions by using an Input-Convex Neural Network~(ICNN)~\citep{amos2017input}. 
    This is a popular architecture, but not the only one ensuring convexity~\citep{gagneux2025convexity}. 
    ICNNS have been used in many studies where a convex approximation is sufficient~\citep{chen2020icnn,rosemberg2024learning,liu2026surrogates,liu2026icnn} 
    because their optimization can be done more efficiently~\citep{duchesne2021convex}.
    \item We can use convexification techniques to work more efficiently with activation functions that are not piecewise-linear, 
    which entail Mixed-Integer Non-Linear Programming (MINLP) embedding formulations.    
    For example, this is possible in the MAiNGO solver~\citep{bongartz2018maingo} through the MeLOn framework~\citep{schweidtmann2020melon}. 
    Convex envelopes for S-shaped activation functions are studied in more detail by~\cite{carrasco2026tightening}.
    \item We can avoid the complexities associated with embedding the neural network in the optimization model by using a weaker form of integration,  
    such as by evaluating the neural network from a callback~\citep{parker2024formulations}, 
    which nevertheless may further limit claims of optimality 
    or limit access to network parameters or gradients by the optimization algorithm.
    \item We can reduce the influence of incorrect predictions of the neural network in the solutions produced by the optimization model in some ways. First, we can use multiple neural networks (an \emph{ensemble}) to map inputs to the average output across those networks~\citep{wang2023ensemble}. Second, we can restrict the inputs and output that can be produced by the model to what is observed in the dataset, such as by only allowing inputs within the convex hull of samples~\citep{shi2024trustregions} or by similarly accounting for both inputs and output~\citep{zhu2026validitydomain}. 
    Finally, 
    in the case that access to the ground truth function is available, 
    we may also generate additional samples to retrain the neural network~\citep{liu2025simulator}. 
\end{enumerate}

In a related approach known as \emph{Decision-Focused Learning}~(DFL),  
predictors such as neural networks 
are trained with the primary goal of improving decision quality of the optimization model rather than predictive accuracy alone~\citep{wilder2019melding,elmachtoub2022smart,mandi2024decision,sadana2025survey},  
which is a form of \emph{end-to-end learning}~\citep{donti2017task}. 
In the context of embedding neural networks in optimization models, this modification may be beneficial, 
provided that repeatedly solving the optimization model during training is sufficiently inexpensive.

There are other forms in which prediction and optimization can be combined. 
First, the optimization model can be represented as a differentiable layer of the neural network~\citep{amos2017optnet,agrawal2019differentiable,ferber2020mipaal}. 
That brings along the cost of solving an optimization problem in inference (forward pass) 
and performing something similar to sensitivity analysis by differentiating on an optimal solution in training (backward pass).
Second, the neural network may be used to produce decisions directly, 
in what is often denoted an \emph{optimization proxy}~\citep{chen2024endtoend,chen2024compact}. 
While it is only possible to ensure a worst-case optimality gap or feasibility straight from the neural network in some special cases, 
feasible solutions can be obtained by postprocessing and optimization proxies 
can be used in many settings where a short runtime is critically important. 

More broadly, machine learning and AI in general have been applied in many different ways to mathematical optimization. 
Some of this work taps into the long-standing application of machine learning to improve optimization solvers and algorithms in general, 
recently discussed in a survey by~\cite{bengio2021machine}, 
but dating back to earlier work such as \cite{hutter2009paramils} and \cite{kadioglu2010isac}. 
In the last years, 
we have seen a growing use of Large Language Models (LLMs) in optimization. 
To mention a few noteworthy cases, 
LLMs have been used to explain why optimization models are infeasible~\citep{chen2024diagnosing}; 
generate models based on problem descriptions~\citep{ahmaditeshnizi2024optimus}; 
configure solver parameters based on problem structure~\citep{lawless2025llms}; 
and explain a model and answer what-if questions to end users~\citep{chen2026optichat}. 

Last, but not least, 
optimization is usually how training in machine learning is framed. 
To mention but one reference for the interested reader, consider \cite{gambella2021optimization}.

\section{Conclusion}

In this tutorial, 
we have shown how neural networks can be beneficial when used for constraint learning, 
which consists of embedding a machine learning model as part of an optimization model. 
We have discussed the fundamentals of this type of approach, 
what technical choices are often made, 
which applications have been explored, 
how to calibrate a neural network for more effective and efficient results, 
and what references to consider for advanced topics related to constraint learning. 
This is currently an active topic of research, 
for which reason we recommend the reader to stay alert to new developments. 
May your problems have good approximations and fast solutions! 

\section*{Acknowledgements}
We would like to thank for the feedback of those who attended presentations of this work the Analytics+ 2025, ICS 2026, and CORS 2026 conferences and the MIP 2026 workshop; as well as in seminars at the University of Iowa's Tippie College of Business and Western University's Ivey Business School.


\beginFEappendix

\section*{Appendix} The sections below provide more references, context, and detailed mathematical formulations. 
They are mentioned when relevant in the main text of the tutorial. 
We believe that they are not necessary in a first read, but that they may become more relevant as the reader keeps refining their understanding and needs further background on the topics.

\section{About CP}\label{ap:about_cp}

Constraint Programming~(CP) is a discipline of modeling and solving combinatorial problems based on backtracking, logical inference, and specialized domain propagation algorithms~\citep{dechter2003constraint}. CP has many parallels with mathematical optimization~\citep{hooker2012integrated}; and tends to be a competitive alternative when obtaining feasible solutions is difficult, such as in scheduling problems~\citep{laborie2018cpoptimizer}. \cite{lustig2001program} present a qualitative discussion of the merits and differences between CP and mathematical optimization. 
While this tutorial is focused on mathematical optimization, we will refer the reader to some papers using CP due to its complementary strengths.

\section{Adopting the Term Constraint Learning in Mathematical Optimization}\label{ap:adopting_cl}

The term \emph{constraint learning} has been widely adopted since 2023 by publications describing mathematical optimization models embedding machine learning models, which in most cases are neural networks, including \cite{alcantara2023chance,alcantara2023distributional,fajemisin2024ocl,sang2024carbonflow,song2024cldispatch,shi2024trustregions,cacciola2024structured,maragno2025miocl,hooshmand2025integration,zhu2026validitydomain,pham2026sparse}, and \cite{tong2026gradient}.

\section{Neural Network Structure and Mapping}\label{ap:nn_structure_map}

The neural network in Figure~\ref{fig:nn_figure} maps every input $\vx = [x_1 ~ x_2 ~ \ldots ~ x_{n_0}]^T \in \sR^{n_0}$ to an output $y \in \sR$. This mapping is supported by $L$ hidden layers, denoted by $\ell = 1$ to $\ell=L$. 
The $\ell$-th hidden layer has $n_{\ell}$ neurons, 
denoted by $i=1$ to $i=n_{\ell}$.  
Those neurons produce intermediary outputs $\vh^\ell = [\vh^\ell_1 ~ \vh^\ell_2 ~ \ldots ~ \vh^\ell_{n_\ell}]^T$. 
For simplicity, we can think of $\vx$ as $\vh^0$. 
Every hidden layer has a matrix of weights $\mW^\ell$ and a bias vector $\vb^\ell$. 
Using the $i$-th column of weights $\mW^\ell_{i *}$ and the $i$-th bias element $b^\ell_i$, 
we calculate the preactivation output $g^\ell_i = \sum\limits_{j=1}^{n_{\ell-1}} \mW^\ell_{i j} h^\ell_j + b^\ell_i$. 
Every hidden layer also has an activation function $\sigma^\ell : \sR \rightarrow \sR$, 
with which the obtain intermediary output $h^\ell_i = \sigma^\ell(g^\ell_i)$.
For the output layer, we similarly use 
a matrix of weights $\mW^{L+1}$ and a bias vector $\vb^{L+1}$, 
which become a vector and scalar since $y$ is a scalar. 
We do not use an activation function for the output; 
and thus, for simplicity, we can think of $y$ as $g^{L+1}_1$.

\section{Why Use a Nonlinear Activation Function?}
\label{app:nonlinear_act}

Every neuron $i$ in layer $\ell$ calculates a preactivation output $g^\ell_i$ based on a linear function on the output $\vh^{\ell-1}$ from the previous layer using the set of weights $W_{i *}^\ell$ and the bias $b_i^\ell$; 
then its output is given by $h^\ell_i = \sigma^\ell(g^\ell_i)$. 
In other words, 
the output of every neuron is the composition of the activation function with a linear transformation of its inputs. 

If we were to bypass the activation function by directly outputting the preactivation output ($h^\ell_i = g^\ell_i$), then the entire neural network would be equivalent to a linear transformation. 
The same would happen with any linear activation function $\sigma^\ell(u) = a u + b$ 
because the composition of two linear functions is also a linear function. 
Hence, we need a nonlinear activation function to make the neural networks more expressive than a linear regression.

\section{On Calculating (Pre) Activation Bounds for Neurons with ReLU Activation}\label{app:relu_bounds}

We would like to find small but valid values for $\bar{H}^\ell_i$ and $H^\ell_i$, 
which should be such that $-\bar{H}^\ell_i \leq g_i^\ell$ and $g_i^\ell \leq H^\ell_i$. 
These are the values commonly referred in the literature as \emph{activation bounds}, 
which we prefer to define as \emph{pre}-activation bounds. 
However, they can only be defined if we bound the inputs of the neural network. 
To simplify notation,  
let $-\bar{H}^0_i$ and $H^0_i$ be finite lower and upper bounds of each input $x_i = h_i^0$ of the neural network.

First we explain how to find the smallest possible values. Because this is going to be as expensive as solving the optimization model in which they depend, we will discuss some alternatives afterwards. 
The smallest possible values can be obtained by solving 
\begin{subequations}
\begin{align}
H^{\ell'}_i = & \max g_i^{\ell'} \\
& \begin{array}{r} (\ref{eq:W}) \\ ~ \end{array} & ~~ \begin{array}{ll} \forall & \ell \in \{1, \ldots, \ell'\}, \\ & i \in \{1, \ldots, n_\ell\} \end{array}  \\
& \begin{array}{r} (\ref{eq:g})\text{ to }(\ref{eq:bin}) \\ ~ \end{array} & ~~ \begin{array}{ll} \forall & \ell \in \{1, \ldots, \ell-1\}, \\ & i \in \{1, \ldots, n_\ell\} \end{array} \\
& -\bar{H}^0_i \leq h_0^i \leq H^0_i & \forall i \in \{1, \ldots, n_0\} 
\end{align}
\label{eq:h_bound}
\end{subequations}
and
\begin{subequations}
\begin{align}
\bar{H}^{\ell'}_i = & - \min g_i^{\ell'} \\
& \begin{array}{r} (\ref{eq:W}) \\ ~ \end{array} & ~~ \begin{array}{ll} \forall & \ell \in \{1, \ldots, \ell'\}, \\ & i \in \{1, \ldots, n_\ell\} \end{array}  \\
& \begin{array}{r} (\ref{eq:g})\text{ to }(\ref{eq:bin}) \\ ~ \end{array} & ~~ \begin{array}{ll} \forall & \ell \in \{1, \ldots, \ell-1\}, \\ & i \in \{1, \ldots, n_\ell\} \end{array} \\
& -\bar{H}^0_i \leq h_0^i \leq H^0_i & \forall i \in \{1, \ldots, n_0\} 
\end{align}
\label{eq:hbar_bound}
\end{subequations}
We assume that $\bar{H}^\ell_i > 0$, since otherwise the neuron is never truly inactive; 
and that $H^\ell_i > 0$, since otherwise the neuron is never truly active. 
In case either assumption fails, 
we can simplify the MILP formulation by fixing the neuron state. 
Note also that we need to solve those models for the first $\ell$ layers before calculating activation bounds for layer $\ell+1$.

Moving to the other extreme, 
one of the simplest ways to obtain valid bound is by calculating the maximum contributions for each item within the summation the defines the preactivation output, through what is commonly denoted as \emph{interval arithmetic}. 
For a neuron $i$ in layer $\ell > 2$, we calculate the bounds as
\begin{align}
H_i^\ell = \sum_{j=1}^{n_{\ell-1}} \max\{ W_{i j}^\ell H_j^{\ell-1}, 0\} + b_i^\ell
\end{align}
\begin{align}
\bar{H}_i^\ell = \sum_{j=1}^{n_{\ell-1}} \min\{ W_{i j}^\ell H_j^{\ell-1}, 0\} + b_i^\ell
\end{align}
For the first layer, given that zero is not necessarily a valid lower bound of each input to the neuron, we use
\begin{align}
H_i^\ell = \sum_{j=1}^{n_{\ell-1}} \max\{ W_{i j}^\ell H_j^{\ell-1}, W_{i j}^\ell \bar{H}_j^{\ell-1}\} + b_i^\ell
\end{align}
\begin{align}
\bar{H}_i^\ell = \sum_{j=1}^{n_{\ell-1}} \min\{ W_{i j}^\ell H_j^{\ell-1}, W_{i j}^\ell \bar{H}_j^{\ell-1}\} + b_i^\ell
\end{align}

Many authors have explored calibrations balancing the tightening of the bounds with the runtime to find such bounds, sometimes aiming for a compromise that would reduce the runtime of solving the optimization model embedding the entire neural network~\citep{cheng2017maximum,serra2018bounding,fischetti2018dnnmilp,grimstad2019surrogate,liu2021algorithms,badilla2023tradeoff,zhao2024horizon,hojny2024message,sosnin2024control}.

\section{Using the Softmax Function in Linear Optimization Models}\label{app:softmax}

The softmax is an activation function used for shaping the output of a neural network as a probability distribution~\citep{bridle1990probabilistic} by taking all inputs of the layer to which it is applied into account to produce the same number of outputs with \[
\sigma^\ell(\vu)_i = \dfrac{e^{u_i}}{\sum\limits_{j=1}^{n_\ell} e^{u_j}}.
\]
Note that this moves away from our premise in the tutorial of working with neural networks having a scalar output. 

While we cannot directly embed an entire neural network with softmax into a linear optimization model, we may still apply linear optimization in some cases where softmax is only used at the last layer. 
For example, this is typically what happens when a neural network with ReLU in all other layers is trained for classification and has an output element for each class. 
One application of linear optimization in that case is network verification against adversarial inputs: 
given an input $\bar{\vx}$ that is correctly classified as belonging to class $i$ (and thus $f(\bar{x}_i) \geq f(\bar{x}_i)$ for every output element $j \neq i$), 
we may look for a small perturbation $\boldsymbol{\varepsilon}$ 
such that $\bar{\vx}+\boldsymbol{\varepsilon}$ gets classified as belonging to another specific class $j$. 
Informally, this relates to solve the optimization problem $\max y_j - y_i$ subject to $\vx = \bar{\vx} + \boldsymbol{\varepsilon}$ and $\vy = f(\vx)$, 
where $y_j > y_i$ (a positive objective function value) implies that the incorrect class $j$ has a higher value than the correct class $i$. 
Because the largest inputs of the softmax layer correspond to the largest outputs of the softmax layer, 
we can drop the softmax layer from the verification problem and replace $y_j$ and $y_i$ with the corresponding inputs of the softmax layer, 
say $g^L_i$ and $g^L_j$, 
and still ensure that an adversarial perturbation will be obtained from solving it to optimality if such a perturbation exists.

\section{Other Formulations for the ReLU}\label{ap:formulations}

At a macroscopic level, 
size is not everything, 
but we need to balance the size of a formulation with the strength of its linear relaxation. 
For example, 
let us consider classic Traveling Salesperson Problem~(TSP), 
where the number of constraints of the Miller--Tucker--Zemlin~(MTZ) formulation is quadratic on the number of cities~\citep{miller1960mtz} 
while the Dantzig--Fulkerson--Johnson~(DFG) formulation has an exponential number~\citep{datzig1954tsp}, 
but the latter can be solved reasonably faster 
because the former has a very weak linear relaxation~\citep{balas1985scheduling}---meaning that it allows many more fractional solutions. 
Drawing an analogy with the TSP, 
formulation~(\ref{eq:W})--(\ref{eq:bin}) would be the equivalent of MTZ, 
but we do not have a practical equivalent of DFG. 
Nevertheless, 
we do know how to use general-purpose lifting techniques to produce a stronger extended formulation, 
as illustrated in \cite{huchette2018phd}, 
and then applied for generating cutting planes in \cite{anderson2020strong} 
and manageably larger, but moderately stronger, formulations in \cite{tsay2021partition} that are generalized in \cite{kronqvist2025psplit}. 

\section{How Does Geometry Constrain the Number of Linear Regions?}\label{ap:geometry}

The neural network in Figure~\ref{fig:pwa}, with two inputs ($n_0=2$) and five neurons in a single hidden layer ($L=1$ and $n_1=5$), 
defines 12 linear regions. 
Given that each linear region is defined by a combination of which neurons are active or not, 
then we could potentially have as many as $2^5=32$ regions for $n_0=2$, $L=1$, and $n_1=5$. 
However, we cannot go further than 16 in this case.

\begin{figure}
    \centering
    \includegraphics[width=0.42\linewidth]{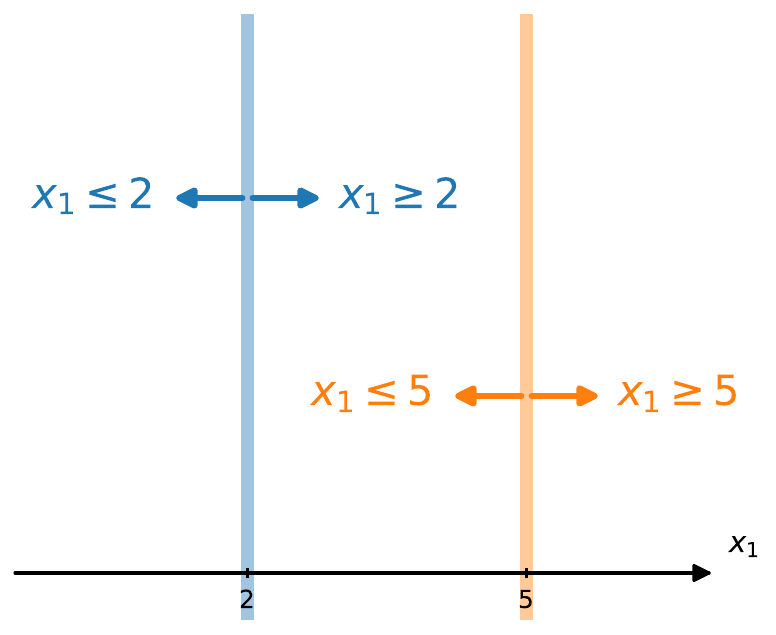}
    \qquad
    \includegraphics[width=0.45\linewidth]{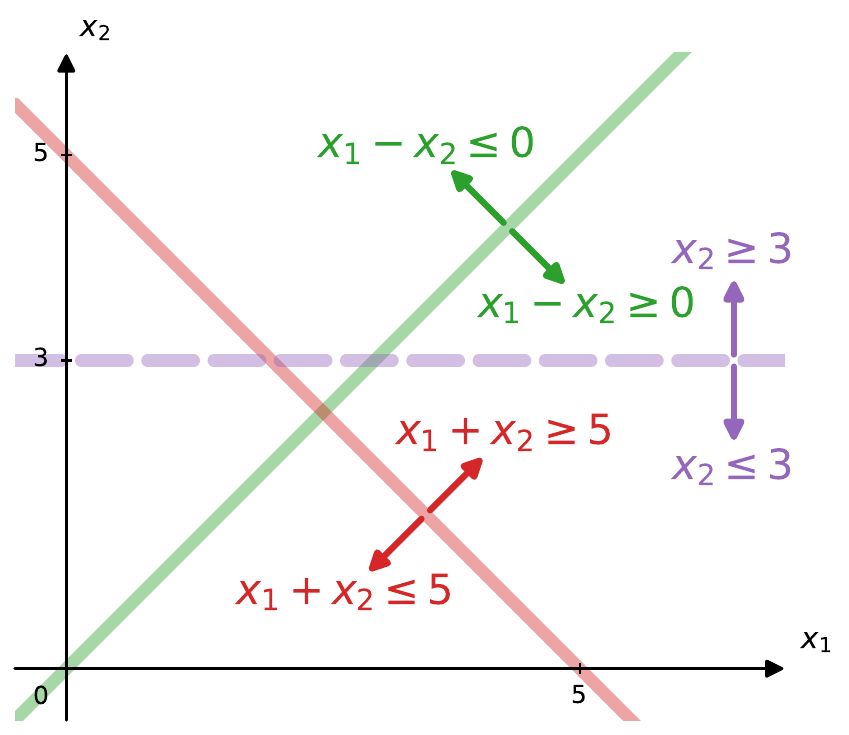}
    \caption{Geometric illustration of how defining half-spaces with hyperplanes in more ways than the number of dimensions does not imply that every intersection of half-spaces is feasible. Left: two partitionings in one dimension lead to a maximum of 3 regions instead of 4. Right: three partitionings in two dimensions lead to a maximum of 7 regions instead of 8.}
    \label{fig:arrangements}
\end{figure}

So why 16 instead of 32?
There is a geometric argument for it, 
which is easier to understand with a simpler example first. 
Consider Figure~\ref{fig:arrangements}, 
where we define half-spaces in one dimension (left) and two dimensions (right). 
In one dimension, 
we can define half-spaces overlapping at the point $x_1=2$ ($x_1 \leq 2$ and $x_1 \geq 2$) 
and at the point $x_1=5$ ($x_1 \leq 5$ and $x_1 \geq 5$). 
However, there is one combination of these half-spaces that is not feasible ($x_1 \leq 2$ and $x_1 \geq 5$), 
for which reason we only define three regions with overlapping boundaries ($x_1 \leq 2$, $2 \leq x_1 \leq 5$, and $x_1 \geq 5$). 
In two dimensions, 
we can now define half-spaces overlapping at lines. 
If we pick two of the lines shown, say $x_1 - x_2 = 0$ and $x_1 + x_2 = 5$, we can see four regions with overlapping boundaries. 
But if we consider the three lines, hence including $x_2=3$, we can see only seven regions instead of eight because one combination is infeasible 
($x_1 - x_2 \geq 0$ and $x_1 + x_2 \leq 5$ with $x_2 \geq 3$). 
Geometrically, 
this is because the third line cannot go through the interior of more than three regions defined by the other two lines. 
More generally, 
if we have a $d$-dimensional space and n equalities defining half-spaces in it, 
the maximum number of regions that we can define is $\sum_{i=0}^{d} \binom{n}{i}$~\citep{zaslavsky1975book}. 
In other words, the dimension of the space affects how many full-dimensional regions we can get by splitting it into half-spaces: 
we can only get $2^n$ regions with $n$ inequalities if $d \geq n$.
Returning to the previous example: 
with $n_0=2$, 
every neuron in layer $1$ draws a line $\left(\sum_{j=1}^{n_1} \mW^1_{i j} x_j + b^1_i = 0\right)$ between the half-space in $x$ where the neuron is active $\left(\sum_{j=1}^{n_1} \mW^1_{i j} x_j + b^1_i \geq 0\right)$ and the half-space where the neuron is inactive $\left(\sum_{j=1}^{n_1} \mW^1_{i j} x_j + b^1_i \leq 0\right)$; 
and with $n_1=5$ neurons, 
we can only define up to $\sum_{i=0}^{n_0} \binom{n_1}{i} = 16$ regions.

\section{Details About Dataset Generation}
\label{app:dataset_generation}

We have generated 10,000 observations, which were split at a 80:20 ratio between the training and validation sets. 

The observation values were obtained as follows. 
For the spending values of the eight advertising channels, 
the values were generated from a uniform distribution with minimum bound 0 and maximum bounds of
120 for television, 
100 for radio, 
300 for social media, 
200 for influencer sponsorship, 
80 for billboard advertising, 
100 for browser advertising, 
40 for email advertising, and 
60 for event sponsorship. 
For the choice of social media, 
we have used 5 binary columns, 
randomly setting all but one to zero (a \emph{one-hot encoding}). 
For the exogenous parameters, 
product price was generated from a normal distribution with mean 500 and standard deviation 20, rounded to two decimal places; 
sales in the last month was generated from a normal distribution with mean 2500 and standard deviation 200; 
and marketing budget of the largest competitor was generate from a uniform distribution with bounds 0 and 1,000. 
The product sales are evaluated using the ground truth function $f$, to be discussed next, to which we added a random noise generated from a normal distribution with mean 0 and standard deviation 1, 
and then we rounded to the nearest integer.

The ground truth function $f$ was defined as follows. 
Let $\vx$ denote the vector of spending values, 
with $x_1$ for television, 
$x_2$ for radio, 
$x_3$ for social media, 
$x_4$ for influencer sponsorship, 
$x_5$ for billboard advertising, 
$x_6$ for browser advertising, 
$x_7$ for email advertising, and 
$x_8$ for event sponsorship. 
Let $\vz \in [0,1]^5$ denote the vector of social media choice. 
Let $\boldsymbol{\xi}$ denote the vector of exogenous parameters, 
with $\xi_1$ for product price, 
$\xi_2$ for sales in the last month, and 
$\xi_3$ for marketing budget of the largest competitor. 
We define the advertising-response component as 
\begin{align*}
R(\vx,\vz) = & 100 \ell(x_1) + 50 \ell(x_2) + 100 \sum_{i=1}^5 z_i \ell\big(g_i(x_3)\big) \\ & + 130 \ell(x_4) + 30 \ell(x_5) + 60 \ell (x_6) \\ & + 90 \ell(x_7) + 20 \ell(x_8)
\end{align*}
where
\[
\ell(u) = \log_2\left(\max\{u,1\}\right)
\]
and 
\[
g_i(u) = a_i u - b_i u^2
\]
with 
$a = [6.66, 2.22, 1.33, 0.95, 0.66]^T$ and $b = [0.111, 0.012, 0.004, 0.002, 0.001]^T$. 
We define the contextual contribution as 
\[
C(\boldsymbol{\xi}) = \frac{10000}{\xi_1} - 0.12 \xi_3
\]
and combine those components into the ground truth function 
\[
f(\vx,\vz, \boldsymbol{\xi}) = 0.2 \xi_2 + 0.8 [ R(\vx,\vz) + C(\boldsymbol{\xi})].
\]

\section{Additional Details About Computational Experiments}\label{ap:experiments}

All numerical experiments were implemented in Python 3.10.8 
using PyTorch 2.13.0 for training the neural networks and 
Gurobi 9.5.2 for solving the MILP models. 
The experiments were conducted on a dual-socket compute node
at the University of Iowa's Argon cluster, 
with 56 physical
CPU cores, 112 hardware threads, and approximately 512 GB
of system memory.
The source code and example dataset are publicly available at
\url{https://github.com/fangjiaxiao/constraint-learning-nn-experiments}.

Each MILP was solved using a single Gurobi thread, with up to
10 independent optimization tasks executed concurrently.
Each MILP was solved with a time limit of 1,000 seconds and a relative MIP gap tolerance of $1\%$. We measured wall-clock time for each optimization task, including loading, MILP construction, and solver execution. 
The maximum budget $B$ was set to 500 and the exogenous parameters were set to their mean value on the dataset.

Each model was trained for up to 1,000 epochs using the Adam optimizer, with a batch size of 64, a learning rate of $10^{-3}$, and a weight decay of $10^{-4}$.
Training minimized the mean squared error, with early stopping after 30 consecutive epochs without improvement in validation MSE. After pruning, each model was fine-tuned for up to 200 epochs using the same optimizer settings, with early stopping after 20 consecutive epochs without improvement in validation MSE. The checkpoint with the lowest validation MSE was retained at each stage. 

\section{Additional Results}\label{ap:results}

\begin{figure}[h]
    \centering
    \includegraphics[width=0.7\linewidth]{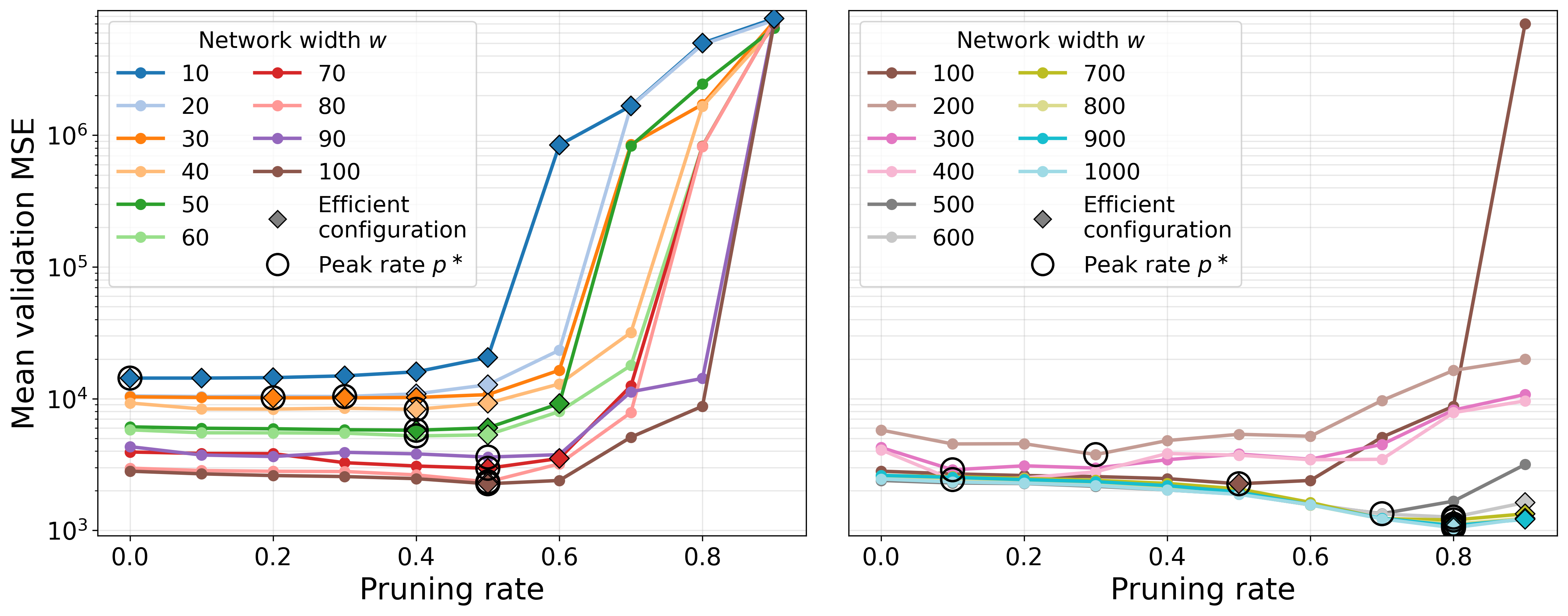}
    \caption{Comparison of predictive quality with pruning rate $p$ (as a line chart for each width $w$) in small (left) and large (right) neural networks. The pruning rate $p*$ with smallest validation error for each width $w$ is circled. The configurations at the efficient frontier are diamond-shaped.}
    \label{fig:pruning_vs_mse}
\end{figure}

Figure~\ref{fig:pruning_vs_mse} complements Figure~\ref{fig:pruning_vs_mse_selected} by showing the results for all network widths used in the experiments. 


\finishFEappendix

\bibliographystyle{abbrvnat}
\bibliography{references} 

\end{document}